\documentclass[11pt,a4paper]{article}
\usepackage[latin1]{inputenc}
\usepackage{amsmath}
\usepackage{amsthm}
\usepackage{amsfonts}
\usepackage{amsfonts,amsthm,latexsym,amsmath,amssymb,amscd,epsfig,psfrag,enumerate}
\usepackage{graphics,graphicx, bezier, float, color, hyperref}
\usepackage{amssymb,url}
\usepackage{multienum}
\usepackage[table]{xcolor}
\usepackage{multicol,multirow}
\usepackage{graphicx}
\usepackage{fancyvrb}
\usepackage{parskip}
\usepackage[toc,page]{appendix}
\usepackage{mathrsfs}
\usepackage[top=2.7cm, bottom=2.7cm, left=1.5cm, right=1.5cm]{geometry}
\usepackage{xcolor}

\makeatletter
\newcommand*{\house}[1]{%
	\mathord{%
		\mathpalette\@house{#1}%
	}%
}
\newcommand*{\@house}[2]{%
	\dimen@=\fontdimen8 %
	\ifx#1\scriptscriptstyle\scriptscriptfont
	\else\ifx#1\scriptstyle\scriptfont
	\else\textfont\fi\fi
	3 %
	\sbox0{%
		$#1%
		\vrule width\dimen@\relax
		\overline{%
			\kern2\dimen@
			\begingroup % to keep changes of \dimen@ in #2 local
			#2%
			\endgroup
			\kern2\dimen@
		}%
		\vrule width\dimen@\relax
		\mathsurround=1.5\dimen@ % outside margin
		$%
	}%
	\ht0=\dimexpr\ht0-\dimen@\relax
	\dp0=\dimexpr\dp0+2\dimen@\relax
	\vbox{%
		\kern\dimen@ % reinsert previously removed space
		\copy0 %
	}%
}

\usepackage{blkarray}
\newtheorem{theorem}{Theorem}[section]
\newtheorem{prop}{Proposition}[section]
\newtheorem{lemma}{Lemma}[section]

\usepackage[none]{hyphenat}[section]
\newtheorem{definition}{Definition}[section]

\numberwithin{equation}{section}
\numberwithin{table}{section}
\numberwithin{figure}{section}

\title{Solutions to a Pillai--type equation involving Tribonacci numbers and $S-$units}
\author{Herbert Batte$^{1,*} $ and Florian Luca$^{2}$}
\date{}

\begin{document}
\maketitle
\abstract{ Let $ \{T_n\}_{n\geq 0} $ be the sequence of Tribonacci numbers. In this paper, we study the exponential Diophantine equation $T_n-2^x3^y=c$, for $n,x,y\in \mathbb{Z}_{\ge0}$. In particular, we show that there is no integer $c$ with at least six representations of the form $T_n-2^x3^y$.} 

{\bf Keywords and phrases}: Tribonacci numbers; $S-$units; linear forms in logarithms; $p-$adic numbers; Pillai's problem.
 
{\bf 2020 Mathematics Subject Classification}: 11B39, 11D61, 11D45, 11Y50.

\thanks{$ ^{*} $ Corresponding author}

\section{Introduction}\label{intro}
\subsection{Background}
\label{sec:1.1}
We consider the sequence of Tribonacci numbers $\{T_n\}_{n\ge 0}$ defined by $T_0=0$, $T_1=T_2=1$ and the recurrence relation $$T_{n+3}=T_{n+2}+T_{n+1}+T_{n},$$ for all $n \geq 0$.
The first few terms of this sequence are given by 
$$
0,\;1,\;1,\;2,\;4,\;7,\;13,\;24,\;44,\;81,\;149,\;274,\;504,\ldots. 
$$
The Diophantine equation
\begin{equation}\label{1.1t}
	a^{x} - b^{y} = c,
\end{equation}
where fixed integers \( a > 1 \), \( b > 1 \) and \( c \) are involved, is known as the Pillai equation. Originating from Pillai's work in \cite{PIL1}, this equation's potential to produce multiple nonnegative integer solutions \( (x, y) \) has been a subject of interest. Pillai's significant contribution was the demonstration that, for positive, coprime integers \( a \) and \( b \), and when \( |c| > c_0(a, b) \), there is no more than at most a single solution \( (x, y) \).

Further explorations of the Pillai problem have focused on the cases where \( a \) is fixed as \( 2 \), \( 3 \), or a particular prime \( p \), replacing the sequence of powers of powers of $b$ with other exponentially growing sequences of positive integers, such as Pell numbers, Fibonacci numbers, Tribonacci numbers, and more complex \( k \)--generalized Fibonacci numbers for an integer  parameter \( k\ge 2 \). These studies, see for example \cite{BAT}, \cite{BFY}, \cite{CPZ}, \cite{DDA}, and \cite{DFR}, have largely confirmed the conclusion of the original Pillai problem: all integers are uniquely represented in the above way, except for finitely many outliers.

In a number field ${\mathbb K}$ with ring of integers ${\mathcal O}_{\mathbb K}$ and a finite set of prime ideals \( S \), an element \( x \in {\mathbb K} \) is an \( S \)--unit if its principal fractional ideal is a product of primes in \( S \). For rational numbers, an \( S \)--unit is a rational number whose numerator and denominator are divisible only by primes in \( S \). More research on the variation of \eqref{1.1t} are found in  \cite{BAT2} and \cite{VZ}. In \cite{VZ}, in equation \eqref{1.1t} Fibonacci numbers replaced the powers of \( a \), and \( S \)--units replaced the powers of \( b \), while in the similar study \cite{BAT2}, equation \eqref{1.1t} was investigated with Lucas numbers and $S$--units. In this note, we revisit \eqref{1.1t} instead with Tribonacci numbers while retaining the \( S \)--units with \( S = \{2, 3\} \).

Therefore, we investigate the exponential Diophantine equation 
\begin{align}\label{1.2t}
	T_n - 2^x3^y = c,
\end{align}
for \( n, x, y \in \mathbb{Z}_{\ge 0} \).

We discard the situation $n = 1$ since
$T_1 = T_2 = 1$. Thus, we always assume that $n=0$, $n \ge 2$. The main results of this paper are the following theorems. 
\subsection{Main Results}\label{sec:1.2t}
\begin{theorem}\label{1.2at} 
	The Diophantine equation \eqref{1.2t} has in the case that $c = 0$, exactly five
	solutions, namely
	$$(n,x,y)= (2, 0, 0), (3, 1, 0), (4, 2, 0),	(7, 3, 1),
	(9, 0, 4). $$
	Furthermore, this representation is given by
	\begin{align*}
		0&=T_2-2^03^0=T_3-2^1 3^0=T_4-2^23^0=T_7-2^3 3^1=T_9-2^0 3^4.
	\end{align*}
\end{theorem}

\begin{theorem}\label{1.2bt} 
	Let $c\in \mathbb{N}$ such that the Diophantine equation \eqref{1.2t} has at least four solutions $(n,x,y)\in \mathbb{Z}_{\ge 0}^3$. Then, $$c=1.$$
	Furthermore, this representation is given by
	\begin{align*}
	1&=T_3-2^03^0=T_4-2^0 3^1=T_5-2^13^1=T_6-2^2 3^1.
	\end{align*}
\end{theorem}

\begin{theorem}\label{1.2ct} 
	Let $c\in -\mathbb{N}$ such that the Diophantine equation \eqref{1.2t} has at least five solutions $(n,x,y)\in \mathbb{Z}_{\ge 0}^3$. Then, $$c \in \{-8,-2\}.$$
	Furthermore, these representations are given by
	\begin{align*}
		-8&=T_0-2^33^0=T_2-2^03^2=T_4-2^2 3^1=T_7-2^53^0=T_{12}-2^9 3^0,\\
		-2&=T_0-2^13^0=T_2-2^03^1=T_3-2^2 3^0=T_4-2^13^1=T_{5}-2^0 3^2.
	\end{align*}
\end{theorem}

\section{Methods}
\subsection{Preliminaries}
The Tribonacci sequence $(T_n)_{n\geq 0}$ has a characteristic polynomial given by
\[
\Psi(X) := X^3 - X^2 - X - 1.
\]
This polynomial $\Psi(X)$ is irreducible in $\mathbb{Q}[X]$, and has a positive real zero
\[
\alpha = \frac{1}{3} \left( 1 + (19 + 3\sqrt{33})^{1/3} + (19 - 3\sqrt{33})^{1/3} \right),
\]
lying strictly outside the unit circle and two complex conjugate zeros $\beta$ and $\gamma$ lying strictly inside the unit circle. Furthermore, $|\beta| = |\gamma| = \alpha^{-1/2}$. According to Dresden and Zu \cite{DZ}, we have
\begin{align}\label{2.1t}
T_n = C_\alpha\alpha^{n-1} + C_\beta\beta^{n-1}  + C_\gamma\gamma^{n-1}  \quad \text{for all } n \geq 0,
\end{align}
where $C_X = (X - 1)/(4X - 6)$. Dresden and Zu \cite{DZ}, also showed that the contribution of the complex conjugate zeros $\beta$ and $\gamma$ to the right--hand side of \eqref{2.1t} is very small. More precisely,
\begin{align}\label{2.2t}
\left|T_n - C_\alpha\alpha^{n-1} \right| < \frac{1}{2} \quad \text{for all } n \geq 0.
\end{align}
The minimal polynomial $P(X)$ of $C_\alpha$ over the integers is given by
\begin{equation}
\label{eq:poly}
P(X)=44X^3 - 44X^2 + 12X - 1,
\end{equation}
and has zeros $C_\alpha$, $C_\beta$, $C_\gamma$ with $|C_\alpha|$, $|C_\beta|$, $|C_\gamma| < 1$. Numerically, these values are approximated as
\begin{align}
\label{eq:approx}
1.83 &< \alpha < 1.84,\nonumber\\
0.73 < |\beta| &= |\gamma| = \alpha^{-1/2} < 0.74,\nonumber\\
0.61 &< |C_\alpha| < 0.62,\\
0.19 < |C_\beta| &= |C_\gamma| < 0.20.\nonumber
\end{align}

Let $\mathbb{K} := \mathbb{Q}(\alpha, \beta)$ be the splitting field of the polynomial $\Psi$ over $\mathbb{Q}$. Then, $[\mathbb{K} : \mathbb{Q}] = 6$. Furthermore, $[\mathbb{Q}(\alpha):\mathbb{Q}]=3$. The Galois group of $\mathbb{K}$ over $\mathbb{Q}$ is given by
\[
\mathcal{G} := \text{Gal}(\mathbb{K}/\mathbb{Q}) \simeq \{ (1), (\alpha\beta), (\alpha\gamma), (\beta\gamma), (\alpha\beta\gamma), (\alpha\gamma\beta) \} \simeq S_3.
\]
Thus, we identify the automorphisms of $\mathcal{G}$ with the permutations of the zeros of the polynomial $\Psi$. For instance, the permutation $(\alpha\gamma)$ corresponds to the automorphism $\sigma : \alpha \mapsto \gamma$, $\gamma \mapsto \alpha$ and $\beta \mapsto \beta$.

We need additional properties. 
\begin{prop}
\label{prop}
The following hold:
\begin{itemize}
\item[\rm{(i)}] There exists a unique $z\in {\mathbb K}$ such that $z^3=\alpha/\beta$.
\item[\rm{(ii)}] If $m\ge 2$, then $\pm z$ is not an $m$th power of some other element in ${\mathbb K}$.
\end{itemize}
\end{prop}

\begin{proof}
(i) For the existence of $z$, see Proposition 5.1 (i) in \cite{BLNOW}. The uniqueness of $z$ follows from the fact that ${\mathbb K}$ does not contain cubic roots of unity. 
The minimal polynomial $Q(X)$ of $z$ is 
$$
Q(X)=X^6+X^5+2X^4+3X^3+2X^2+X+1.
$$
(ii) Assume that $\pm z=y^m$ for some integer $m\ge 2$. We may assume that $m=p$ is a prime and write $y=(\pm z)^{1/p}\in {\mathbb K}$. Note that $y$ is an algebraic integer as it is a unit. 
Also, $|y|=|z|^{1/p}=|\alpha/\beta|^{1/3p}=|\alpha|^{1/2p}$. Note that all conjugates of $y$ have absolute values in the set $\{|\alpha|^{1/2p}, 1, |\alpha|^{-1/2p}\}$. 
Using the main result of Voutier in \cite{Vo}, we have that 
$$
|\alpha|^{1/2p}=|y|	\ge 1+\frac{1}{2d}\left(\dfrac{\log\log d}{\log d}\right)^3,
$$
where $d$ is the degree of ${\mathbb K}$. 
Using $d=6$ and taking logarithms, we get
$$
p\le  \frac{\log\alpha}{2\log\left(1+\dfrac{1}{2\cdot 6}\left(\dfrac{\log\log 6}{\log 6}\right)^3\right)}<106.
$$
Now we used SageMath to check that for prime $p\le 106$, the polynomial $Q(X^p)$ is irreducible. This proves (ii). 
\end{proof}

We recall one additional simple fact from calculus. This is Lemma 1 in \cite{VZ}.
\begin{lemma}[Lemma 1 in \cite{VZ}]\label{lem2.1t}
If $x\in \mathbb{R}$ satisfies $|x|<\frac{1}{2}$, then $|\log(1+x)|<\frac{3}{2}|x|$.	
\end{lemma}

To end this section, we present an analytic argument which is Lemma 7 in \cite{GL}. 
\begin{lemma}[Lemma 7 in \cite{GL}]\label{Guzt} If $ m \geq 1 $, $T > (4m^2)^m$ and $T > \displaystyle \frac{z}{(\log z)^m}$, then $$z < 2^m T (\log T)^m.$$	
\end{lemma}

\subsection{Linear forms in logarithms}
We use several times Baker--type lower bounds for nonzero linear forms in four or more logarithms of algebraic numbers. There are many such bounds mentioned in the literature like that of Baker and W{\"u}stholz from \cite{BW} or Matveev from \cite{MAT}. Before we can formulate such inequalities we need the notion of height of an algebraic number recalled below.

\begin{definition}\label{def2.1t}
	Let $ \lambda $ be an algebraic number of degree $ d $ with minimal primitive polynomial over the integers $$ a_{0}x^{d}+a_{1}x^{d-1}+\cdots+a_{d}=a_{0}\prod_{i=1}^{d}(x-\lambda^{(i)}), $$ where the leading coefficient $ a_{0} $ is positive. Then, the logarithmic height of $ \lambda$ is given by $$ h(\lambda):= \dfrac{1}{d}\Big(\log a_{0}+\sum_{i=1}^{d}\log \max\{|\lambda^{(i)}|,1\} \Big). $$
\end{definition}
 In particular, if $ \lambda$ is a rational number represented as $\lambda:=p/q$ with coprime integers $p$ and $ q\ge 1$, then $ h(\lambda ) = \log \max\{|p|, q\} $. 
The following properties of the logarithmic height function $ h(\cdot) $ will be used in the rest of the paper without further reference:
\begin{equation}\nonumber
	\begin{aligned}
		h(\lambda_{1}\pm\lambda_{2}) &\leq h(\lambda_{1})+h(\lambda_{2})+\log 2;\\
		h(\lambda_{1}\lambda_{2}^{\pm 1} ) &\leq h(\lambda_{1})+h(\lambda_{2});\\
		h(\lambda^{s}) &= |s|h(\lambda)  \quad {\text{\rm valid for}}\quad s\in \mathbb{Z}.
	\end{aligned}
\end{equation}

A linear form in logarithms is an expression
\begin{equation}
	\label{eq:Lambda}
	\Lambda:=b_1\log \lambda_1+\cdots+b_t\log \lambda_t,
\end{equation}
where for us $\lambda_1,\ldots,\lambda_t$ are positive real  algebraic numbers and $b_1,\ldots,b_t$ are nonzero integers. We assume, $\Lambda\ne 0$. We need lower bounds 
for $|\Lambda|$. We write ${\mathbb L}:={\mathbb Q}(\lambda_1,\ldots,\lambda_t)$ and $D$ for the degree of ${\mathbb L}$.
We start with the general form due to Matveev \cite{MAT}. 

\begin{theorem}[Matveev, \cite{MAT}]
	\label{thm:Matt} 
	Put $\Gamma:=\lambda_1^{b_1}\cdots \lambda_t^{b_t}-1=e^{\Lambda}-1$. Assume $\Gamma\ne 0$. Then 
	$$
	\log |\Gamma|>-1.4\cdot 30^{t+3}\cdot t^{4.5} \cdot D^2 (1+\log D)(1+\log B)A_1\cdots A_t,
	$$
	where $B\ge \max\{|b_1|,\ldots,|b_t|\}$ and $A_i\ge \max\{Dh(\lambda_i),|\log \lambda_i|,0.16\}$ for $i=1,\ldots,t$.
\end{theorem}

We also employ a $p$--adic variation of Laurent's result as established by Bugeaud and Laurent in \cite{BL}, Corollary 1. Prior to outlining their result, we require a few additional notations. 

\begin{definition}\label{def2.2t}
Let \( p \) be a prime number. The \( p \)-adic valuation of an integer \( x \), denoted by $\nu_p(x)$,  is defined by
\[
\nu_p(x) := 
\begin{cases} 
	\max\{k \in \mathbb{N} : p^k \mid x\}, & \text{if } x \neq 0; \\
	\infty, & \text{if } x = 0.
\end{cases}
\]
Additionally, if $x=a/b$ is a rational number and $a,~b$ are integers, then we put 
$$\nu_p(x):=\nu_p(a)-\nu_p(b).
$$
\end{definition}

The formula for $\nu_p(x)$ when $x$ is rational given by Definition \ref{def2.2t} does not depend on the representation of $x$ as a ratio of integers $a/b$. It follows easily from Definition \ref{def2.2t} that if $x$ is rational, then
$$
\nu_p(x)=\text{ord}_p(x),
$$
where $\text{ord}_p(x)$ is the exponent of $p$ in the factorization of $x$. For example, $\nu_2(9/8)=-3$. Now for the algebraic number $\lambda$ in Definition \ref{def2.1t}, we define
\[
\nu_p(\lambda):=\frac{\nu_p(a_d/a_0)}{d},
\]
that is, it is the $p$--adic valuation of the rational number $a_d/a_0$ divided by the degree $d$ of $\lambda$. For example, when $x$ is a rational number, we write it as $x=a_d/a_0$ with coprime integers $a_d$ and $a_0\ge 1$, then 
its minimal polynomial is $f(X)=a_0X-a_d$ has degree $1$ so 
$$
\nu_p(x)=\nu_p(a_d/a_0),
$$
consistent with Definition \ref{def2.2t}. 
The \( p \)--adic valuation gives rise to the conventional absolute value. When the rational numbers \( \mathbb{Q} \) are completed with the standard absolute value, the outcome is the set of real numbers, \( \mathbb{R} \). Conversely, employing the \( p \)-adic absolute value for the completion of \( \mathbb{Q} \) yields the \( p \)-adic numbers, represented as \( \mathbb{Q}_p \).

Similarly to the preceding context, let $\lambda_{1}$ and $\lambda_{2}$ be algebraic numbers over $\mathbb{Q}$, treated as elements of the field ${\mathbb K}_p: = \mathbb{Q}_p(\lambda_{1},\lambda_{2})$, with $D := [\mathbb{Q}_p(\lambda_{1},\lambda_{2}):\mathbb{Q}_p]$. Similar to the situation in Theorem \ref{thm:Matt} above, we must employ an adjusted height function. Specifically, we express it as follows:
\[h'(\lambda_{i}) \ge \max \left\{h(\lambda_{i}),\dfrac{\log p}{D}\right\}, ~~\text{for }~i=1,2.\]
\begin{theorem}[Bugeaud and Laurent, \cite{BL}]\label{lem:Bugt}
	Let $b_1$, $b_2$ be positive integers and suppose that $\lambda_{1}$ and $\lambda_{2}$ are multiplicatively independent algebraic numbers such that $\nu_p(\lambda_{1})=\nu_p(\lambda_{2})=0$. Put 
	$$  E':=\dfrac{b_1}{h'(\lambda_{2})}+\dfrac{b_2}{h'(\lambda_{1})},   $$
	and 
	$$ E:=\max\left\{\log E'+\log\log p+0.4, 10, 10\log p\right\}.    $$
	Then
	$$ \nu_p\left(\lambda_{1}^{b_1}\lambda_{2}^{b_2}-1\right)\le \frac{24pg}{(p-1)(\log p)^4} E^2D^4   h'(\lambda_{1})h'(\lambda_{2}  ), $$
	where $g>0$ is the smallest integer such that $\nu_p\left(\lambda_{i}^g-1\right)>0$.
\end{theorem}
Applying Theorem \ref{thm:Matt} and \ref{lem:Bugt}, we get upper bounds on our variables. 

However,  such upper bounds are too large, thus there is need to reduce them. In this paper, we use the following result related with continued fractions (see Theorem 8.2.4 in \cite{ME}).

\begin{lemma}[Legendre]\label{lem:Legt} Let $ \mu $ be an irrational number, $[a_0,a_1,a_2,\ldots]$ be the continued fraction expansion of $\mu$. Let $p_i/q_i=[a_0,a_1,a_2,\ldots,a_i]$, for all $i\ge 0$, be all the convergents of the continued fraction of $ \mu$, and $M$ be a positive integer. Let $ N $ be a non-negative integer such that
	$ q_{N} > M $.
	Then putting $ a(M) := \max \{a_{i}: i=0,1,2,\ldots,N   \}$, the inequality
	$$ \bigg| \mu-\frac{r}{s}   \bigg| > \dfrac{1}{(a(M)+2)s^2},  $$
	holds for all pairs $ (r, s) $ of positive integers with $ 0 < s < M $. 
\end{lemma}
However, since there are no methods based on continued fractions to find a lower bound for linear forms in more than two variables with bounded integer coefficients, we use at some point a method based on the LLL--algorithm. We next explain this method.

\subsection{Reduced Bases for Lattices and LLL--reduction methods}\label{sec2.3}

Let \( k \) be a positive integer. A subset \( \mathcal{L} \) of the \( k \)--dimensional real vector space \( \mathbb{R}^k \) is called a lattice if there exists a basis \( \{b_1, b_2, \ldots, b_k \}\) of \( \mathbb{R}^k \) such that
\begin{align*}
\mathcal{L} = \sum_{i=1}^{k} \mathbb{Z} b_i = \left\{ \sum_{i=1}^{k} r_i b_i \mid r_i \in \mathbb{Z} \right\}.
\end{align*}
In this situation we say that \( b_1, b_2, \ldots, b_k \) form a basis for \( \mathcal{L} \), or that they span \( \mathcal{L} \). We
call \( k \) the rank of \( \mathcal{L} \). The determinant \( \text{det}(\mathcal{L}) \), of \( \mathcal{L} \) is defined by
\begin{align*}
	 \text{det}(\mathcal{L}) = | \det(b_1, b_2, \ldots, b_k) |,
\end{align*}
with the \( b_i \) being written as column vectors. This is a positive real number that does not depend on the choice of the basis (see \cite{Cas} Sect 1.2).

Given linearly independent vectors \( b_1, b_2, \ldots, b_k \) in \( \mathbb{R}^k \), we refer back to the Gram--Schmidt orthogonalization technique. This method allows us to inductively define vectors \( b^*_i \) (with \( 1 \leq i \leq k \)) and real coefficients \( \mu_{i,j} \) (for \( 1 \leq j \leq i \leq k \)). Specifically,
\begin{align*}
	b^*_i &= b_i - \sum_{j=1}^{i-1} \mu_{i,j} b^*_j,~~~
  \mu_{i,j} = \dfrac{\langle b_i, b^*_j\rangle }{\langle b^*_j, b^*_j\rangle},
\end{align*}
where \( \langle \cdot , \cdot \rangle \)  denotes the ordinary inner product on \( \mathbb{R}^k \). Notice that \( b^*_i \) is the orthogonal projection of \( b_i \) on the orthogonal complement of the span of \( b_1, \ldots, b_{i-1} \), and that \( \mathbb{R}b_i \) is orthogonal to the span of \( b^*_1, \ldots, b^*_{i-1} \) for \( 1 \leq i \leq k \). It follows that \( b^*_1, b^*_2, \ldots, b^*_k \) is an orthogonal basis of \( \mathbb{R}^k \). 
\begin{definition}
The basis \( b_1, b_2, \ldots, b_n \) for the lattice \( \mathcal{L} \) is called reduced if
\begin{align*}
	\| \mu_{i,j} \| &\leq \frac{1}{2}, \quad \text{for} \quad 1 \leq j < i \leq n,~~
	\text{and}\\
	\|b^*_{i}+\mu_{i,i-1} b^*_{i-1}\|^2 &\geq \frac{3}{4}\|b^*_{i-1}\|^2, \quad \text{for} \quad 1 < i \leq n,
\end{align*}
where \( \| \cdot \| \) denotes the ordinary Euclidean length. The constant \( \frac{3}{4} \) above is arbitrarily chosen, and may be replaced by any fixed real number \( \mathcal{P} \) with \( \frac{1}{4} < \mathcal{P} < 1 \),\rm{ (see \cite{LLL} Sect 1)}.
\end{definition}\noindent
Let $\mathcal{L}\subseteq\mathbb{R}^k$ be a $k-$dimensional lattice  with reduced basis $b_1,\ldots,b_k$ and denote by $B$ the matrix with columns $b_1,\ldots,b_k$. 
We define
\[
l\left( \mathcal{L},v\right)= \left\{ \begin{array}{c}
	\min_{u\in \mathcal{L}}||u-v|| \quad  ;~~ v\not\in \mathcal{L}\\
\min_{0\ne u\in \mathcal{L}}||u|| \quad  ;~~ v\in \mathcal{L}
\end{array}
\right.,
\]
where $||\cdot||$ denotes the Euclidean norm on $\mathbb{R}^k$. It is well known that, by applying the
LLL-algorithm, it is possible to give in polynomial time a lower bound for $l\left( \mathcal{L},v\right)\ge c_1$ (see \cite{SMA}, Sect. V.4).
\begin{lemma}\label{lem2.5t}
	Let $v\in\mathbb{R}^k$ and $z=B^{-1}v$ with $z=(z_1,\ldots,z_k)^T$. Furthermore, 
	\begin{enumerate}[(i)]
		\item if $v\not \in \mathcal{L}$, let $i_0$ be the largest index such that $z_{i_0}\ne 0$ and put $\sigma:=\{z_{i_0}\}$, where $\{\cdot\}$ denotes the distance to the nearest integer.
		\item if $v\in \mathcal{L}$, put $\sigma:=1$.
	\end{enumerate}
Finally, we have 
\[
c_1:=\max\limits_{1\le j\le k}\left\{\dfrac{||b_1||}{||b_j^*||}\right\}
\qquad\text{and}\qquad
c_2:= c_1^{-1}\sigma||b_1||.
\]
\end{lemma}

In our application, we are given real numbers $\eta_0,\eta_1,\ldots,\eta_k$ which are linearly independent over $\mathbb{Q}$ and two positive constants $c_3$ and $c_4$ such that 
\begin{align}\label{2.9}
|\eta_0+a_1\eta_1+\cdots +a_k \eta_k|\le c_3 \exp(-c_4 H),
\end{align}
where the integers $a_i$ are bounded as $|a_i|\le A_i$ with $A_i$ given upper bounds for $1\le i\le k$. We write $A_0:=\max\limits_{1\le i\le k}\{A_i\}$. 

The basic idea in such a situation, from \cite{Weg}, is to approximate the linear form \eqref{2.9} by an approximation lattice. So, we consider the lattice $\mathcal{L}$ generated by the columns of the matrix
$$ \mathcal{A}=\begin{pmatrix}
	1 & 0 &\ldots& 0 & 0 \\
	0 & 1 &\ldots& 0 & 0 \\
	\vdots & \vdots &\vdots& \vdots & \vdots \\
	0 & 0 &\ldots& 1 & 0 \\
	\lfloor M\eta_1\rfloor & \lfloor M\eta_2\rfloor&\ldots & \lfloor M\eta_{k-1}\rfloor& \lfloor M\eta_{k} \rfloor
\end{pmatrix} ,$$
where $M$ is a large constant usually of the size of about $A_0^k$ . Let us assume that we have an LLL--reduced basis $b_1,\ldots, b_k$ of $\mathcal{L}$ and that we have a lower bound $l\left(\mathcal{L},v\right)\ge c_1$ with $v:=(0,0,\ldots,-\lfloor M\eta_0\rfloor)$. Note that $ c_2$ can be computed by using the results of Lemma \ref{lem2.5t}. Then, with these notations the following result  is Lemma VI.1 in \cite{SMA}.
\begin{lemma}[Lemma VI.1 in \cite{SMA}]\label{lem2.6t}
	Let $S:=\displaystyle\sum_{i=1}^{k-1}A_i^2$ and $T:=\dfrac{1+\sum_{i=1}^{k}A_i}{2}$. If $c_2^2\ge T^2+S$, then inequality \eqref{2.9} implies that we either have $a_1=a_2=\cdots=a_{k-1}=0$ and $a_k=-\dfrac{\lfloor M\eta_0 \rfloor}{\lfloor M\eta_k \rfloor}$, or
	\[
	H\le \dfrac{1}{c_4}\left(\log(Mc_3)-\log\left(\sqrt{c_2^2-S}-T\right)\right).
	\]
\end{lemma}

SageMath 9.5 is used to perform all the computations in this work.

\subsection{Equations in $\alpha,\beta,\gamma$}

In this section, we prove the following results which will be helpful and recalled later.
\begin{lemma}\label{lem:gamt} The equation 
	\begin{align}\label{gamt}
		\dfrac{1-\gamma^{v}}{1-\gamma^{u}} =\dfrac{1-\alpha^{v}}{1-\alpha^{u}}.
	\end{align}
has no integer solutions $u,v$, with $u>v\ge 1$.
\end{lemma}
\begin{proof}
Let $u,v\in \mathbb{Z}$, with $u>v\ge 1$. Then we can write \eqref{gamt} as
	\begin{align*}
		\alpha^{u}-\alpha^{v}&=\gamma^{u}-\gamma^{v}+\gamma^{v}\alpha^{u}-\gamma^{u}\alpha^{v};\\
		\alpha^{u}-\alpha^{v}=\left|\alpha^{u}-\alpha^{v}\right|&=\left|\gamma^{u}-\gamma^{v}+\gamma^{v}\alpha^{u}-\gamma^{u}\alpha^{v}\right|\\
		&\le |\gamma|^{u}+|\gamma|^{v}+\alpha^{u-0.5v}+\alpha^{v-0.5u},
	\end{align*}
	where we have used the fact that $|\gamma| = \alpha^{-1/2} $. Assume for a moment that $v\ge 5$ which implies that $u\ge 6$ since $u>v$. Then the above inequality implies
	\begin{align*}
		\alpha^{u-2}+\alpha^{u-3} &=
		\alpha^{u}-\alpha^{u-1} \le \alpha^{u}-\alpha^{v}\le 1+\alpha^{u-0.5v}+\alpha^{v-0.5u}\\
		&\le  1+\alpha^{u-2.5}+\alpha^{u-1-0.5\cdot6}=
		1+\alpha^{u-2.5}+\alpha^{u-4}\\
		&<1+\alpha^{u-2.5}+\alpha^{u-3}.
	\end{align*}
	Therefore, $\alpha^{u-2}<1+\alpha^{u-2.5}$. From this, we see that $\alpha^{u-2}\left(1-\alpha^{-0.5}\right)<1$. This implies that
	\begin{align*}
		\alpha^{u-2}	&<\dfrac{\sqrt{\alpha}}{\sqrt{\alpha}-1}<\dfrac{\sqrt{1.84}}{\sqrt{1.83}-1}<3.85,
	\end{align*}
	so that $u<4.24$. This contradicts the assumption that $u\ge 6$. Therefore \eqref{gamt} has no solutions with $v\ge 5$. This means that \eqref{gamt} may have solutions only for $1\le v\le 4$. We go back and rewrite \eqref{gamt} as
	\begin{align*}
	\alpha^{u}-\alpha^{v}&=\gamma^{u}-\gamma^{v}+\gamma^{v}\alpha^{u}-\gamma^{u}\alpha^{v};\\
	\alpha^{u}&=\gamma^{u}-\gamma^{v}+\gamma^{v}\alpha^{u}-\gamma^{u}\alpha^{v}+\alpha^{v}\\
	&\le |\gamma|^{u}+|\gamma|^{v}+\alpha^{u-0.5v}+\alpha^{v-0.5u}+\alpha^{v}\\
		&< 2 + \alpha^{u-0.5}+\alpha^{3.5}+\alpha^{4},
\end{align*}	
so that $\alpha^u\left(1-\alpha^{-0.5}\right)<22$, which gives $u\le 7$.	With help of SageMath, we check for $v\in \{1,2,3,4\}$ and $u\in \{2,3,4,5,6,7\}$ that satisfy \eqref{gamt} and find no solutions.
\end{proof}
\begin{lemma}\label{lem:a/b} 
Assume that the algebraic numbers
	\begin{align*}
\lambda_{1}:=\dfrac{\alpha}{\beta}\qquad {\text{\rm and}}\qquad \lambda_{2}:=\dfrac{C_\alpha}{C_\beta}\left(\dfrac{\left(1-\alpha^{n_1-n}\right)\left(1-\gamma^{n_2-n}\right)-\left(1-\alpha^{n_2-n}\right)\left(1-\gamma^{n_1-n}\right)}{\left(1-\beta^{n_1-n}\right)\left(1-\gamma^{n_2-n}\right)-\left(1-\beta^{n_2-n}\right)\left(1-\gamma^{n_1-n}\right)}\right),
\end{align*}
are multiplicatively dependent for $n>n_1>n_2$. Then $\lambda_2=\pm z^t$ for some integer $t$ with $|t|<80n+244$ and $z$ is defined in Proposition \ref{prop}. 
\end{lemma}
\begin{proof}
Suppose $\lambda_{1}$ and $\lambda_{2}$ are multiplicatively dependent. This means that there exists integers $A$ and $B$, not both zero, such that $\lambda_{1}^A=\lambda_{2}^B$. Since $\lambda_1=z^3$, we have 
$$
\lambda_2^B=z^{3A}.
$$
Let $d:=\gcd(3A,B)$ and write $3A=A_1d,~B=B_1d$. Then $\lambda_2^{B_1}=\zeta_d z^{A_1}$ where $\zeta_d$ is a root of unity of order dividing $d$. Since ${\mathbb K}$ does not contain roots of unity other than $\pm 1$, 
it follows that $\lambda_2^{B_1}=\pm z^{A_1}$. We show that $B_1=1$. Indeed, since $\gcd(A_1,B_1)=1$, we get that $A_1p+B_1q=1$ for some integers $p,~q$. Thus, 
$$
z=z^{A_1p+B_1q}=(\pm 1)^p(\lambda_2^pz^q)^{B_1}.
$$
The above calculation shows that $\pm z$ is a power of exponent $B_1$ of some other number in ${\mathbb K}$, so by (ii) of Proposition \ref{prop}, we get $B_1=1$. Thus, $\lambda_2=\pm z^t$ with 
$t=A_1$, and we need to bound $t$. Note that $|t|=|\log |\lambda_2||/|\log z|$. We estimate 
\begin{eqnarray*}
	|\log|\lambda_{2}||&= & \left|\log\left|\dfrac{C_\alpha}{C_\beta}\left(\dfrac{\left(1-\alpha^{n_1-n}\right)\left(1-\gamma^{n_2-n}\right)-\left(1-\alpha^{n_2-n}\right)\left(1-\gamma^{n_1-n}\right)}{\left(1-\beta^{n_1-n}\right)\left(1-\gamma^{n_2-n}\right)-\left(1-\beta^{n_2-n}\right)\left(1-\gamma^{n_1-n}\right)}\right)\right|\right|\\
	&\le &  | \log |C_\alpha|| +|\log |C_\beta| |+\left|\log \left|\left(1-\alpha^{n_1-n}\right)\left(1-\gamma^{n_2-n}\right)-\left(1-\alpha^{n_2-n}\right)\left(1-\gamma^{n_1-n}\right)\right|\right|\\
	&~& +\left|\log \left|\left(1-\beta^{n_1-n}\right)\left(1-\gamma^{n_2-n}\right)-\left(1-\beta^{n_2-n}\right)\left(1-\gamma^{n_1-n}\right)\right|\right|.
\end{eqnarray*}	
The first two terms above are bounded as $|\log |C_{\alpha}||<0.5$ and $|\log|C_{\beta}||<1.67$. The next two terms are  absolute values of logarithms of absolute values of differences of complex numbers. 
These can be very small or can be very large. In the ``large" part it is easy. Since $\alpha>1$ and $n>n_1>n_2$, $|\beta|=|\gamma|=\alpha^{-1/2}$, we have by the absolute value inequality and the fact that $|1-\alpha^{n_i-n}|=1-\alpha^{n_i-n}<1$ for $i=1,2$, that
\begin{eqnarray}
\label{eq:smallupper1}
 && |(1-\alpha^{n_1-n})(1-\gamma^{n_2-n})-(1-\alpha^{n_2-n})(1-\gamma^{n_1-n})|<|1-\gamma^{n_2-n}|+|1-\gamma^{n_1-n}|\nonumber\\
 & \le & 1+\alpha^{(n-n_1)/2}+1+\alpha^{(n-n_2)/2}<2\alpha^{(n-n_2)/2+2}<\alpha^{n+4},
 \end{eqnarray}
 where we used $\alpha^3>\alpha+1$ and $2<\alpha^2$. Furthermore,
 \begin{eqnarray}
 \label{eq:smallupper2}
&&  |(1-\beta^{n_1-n})(1-\gamma^{n_2-n})-(1-\beta^{n_2-n})(1-\gamma^{n_1-n})|\nonumber\\
& < & |1-\beta^{n_1-n}||1-\gamma^{n_2-n}|+|(1-\beta^{n_2-n}||1-\gamma^{n_1-n}|<2(\alpha^{(n-n_1)/2}+1)(\alpha^{(n-n_2)/2}+1)\nonumber\\
& < & 2 (\alpha^{(n-n_2)/2+2})^2=2\alpha^{n-n_2+4}<\alpha^{n+5}.
\end{eqnarray}
The challenge is computing lower bounds. Let 
$$
u:=1-\beta^{n_1-n},\quad v:=1-\gamma^{n_2-n}.
$$
Then we need to compute a lower bound for 
$$
|(uv)^{(\sigma)}-(uv)^{(\tau)}|
$$
over all the distinct conjugates of $uv$. For example, taking $\sigma={\text{\rm id}}$ and $\tau$ to be the map that swaps $\beta$ with $\gamma$, we get a lower bound on the absolute value the expression appearing in 
the denominator of $\lambda_2$ in the statement of Lemma \ref{lem:a/b} (ignoring $C_{\beta}$), while taking $\sigma$ to be the map that swaps $\alpha$ with $\beta$ and $\tau$ the three cycle $\beta\mapsto\gamma\mapsto\alpha$ we get a lower bound on the numerator appearing in the expression of $\lambda_2$ in the statement of Lemma \ref{lem:a/b} (ignoring $C_{\alpha}$). For this, we use a well--known root separation result of Mahler \cite{Du}. In our particular case, it says that if 
$$
R(X)=\prod_{i=1}^d (X-z_i)\in {\mathbb Z}[X]
$$
is a monic polynomial with integer coefficients and distinct roots then writing 
$$
M(R):=\prod_{i=1}^ d \max\{1, |z_i|\}
$$
for its Mahler measure, we have
$$
|z_i-z_j|>\frac{1}{d^{d/2} M(P)^{d-1}}.
$$
For us,
$$
R(X)=\prod_{\sigma\in S_3} (X-(uv)^{(\sigma)})
$$
has degree $d=6$. Writing down the 6 conjugates $(uv)^{(\sigma)}$, we get a total of $12$ factors namely two occurrences of each of $1-\delta^{n_s-n}$ for each $s\in \{1,2\}$ and $\delta\in \{\alpha,\beta,\gamma\}$. 
Ignoring the contributions of $1-\alpha^{n_s-n}$ (as they are in $(0,1)$) and using $|1-\delta^{n_s-n}|\le 1+\alpha^{(n-n_s)/2}<\alpha^{(n-n_s+4)/2}$ as before twice for each of $s\in \{1,2\}$ and $\delta\in \{\beta,\gamma\}$, we get that
$$
M(P)<(\alpha^{(n-n_1+4)/2})^4(\alpha^{(n-n_2+4)/2})^4=\alpha^{4n-2n_2-2n_1+16}\le \alpha^{4n+10}.
$$
Thus,
$$
|(uv)^{(\sigma)}-(uv)^{(\tau)}|>\frac{1}{6^3 (\alpha^{4n+10})^5},
$$
showing that if the left--hand side above is $<1$, then 
$$
|\log|((uv)^{(\sigma)}-(uv)^{(\tau)})||<\log( 6^3 \alpha^{20n+50})<(20 n+50)\log \alpha+3\log 6.
$$
Comparing the above with \eqref{eq:smallupper1} and \eqref{eq:smallupper2}, we conclude that the last bound above holds in all cases. Putting everything together, we get
$$
|\log |\lambda_2||<0.5+1.67+2((20n+50)\log \alpha+3\log 6)<(40n+100)\log \alpha+13.
$$
We thus get
$$
|t|\le \frac{|\log |\lambda_2||}{\log |z|}\le \frac{(40n+100)\log \alpha+13}{0.5\log \alpha}<80n+244.
$$
\end{proof}

\begin{lemma}
\label{zLucas}
For any positive integer $t$, we have 
\begin{itemize}
\item[\rm{(i)}] 
$
\nu_2(z^t-1)<1+\nu _2(n), \quad\text{and}
$
\item[\rm{(ii)}] 
$
\nu_3(z^t- 1)\le 1+\nu_3(n).
$
\end{itemize}
\end{lemma}

\begin{proof}
For a prime ideal $\pi$ of ${\mathcal O}_{\mathbb K}$ and an element $x\in {\mathbb K}$, we write $\nu_{\pi}(x)$ for the exponent at which $\pi$ appears in the factorization of the principal 
fractional ideal $x{\mathcal O}_{\mathbb K}$. We shall use the following well-known fact. Let $p\in {\mathbb Z}$ be the prime such that $\pi$ sits above $p$ (so, $\pi\mid p$) and let $e_{\pi}:=\nu_{\pi}(p)$ is the ramification index of $\pi$. If $\zeta\in {\mathbb K}$ is such that 
\begin{equation}
\label{eq:valuation}
\nu_{\pi}(\zeta-1)\ge \frac{e_{\pi}}{p-1},\qquad {\text{\rm then}}\qquad \nu_{\pi}(\zeta^t-1)=\nu_{\pi}(\zeta-1)+\nu_{\pi}(n)
\end{equation}
(see Lemma 1 in \cite{PosSch} or Lemma 4.4 in \cite{BiluLuca}). 
We start with $p=2$. Then $2{\mathcal O}_{\mathbb K}=\pi^3$, where $\pi$ is a prime ideal of ${\mathcal O}_{\mathbb K}$. In fact, $\pi$ is principal and generated by $1+z+z^2$, as it can be seen by rewriting the equation 
$$
Q(z)=0\qquad {\text{\rm as}}\qquad (z^2+z+1)^3=2z(z+1)^2(z^2+1),
$$
and using that $z+1,~z^2+1$ are units in ${\mathcal O}_{\mathbb K}$. Now ${\mathcal O}_{\mathbb K}/\pi$ is a field with $4$ elements. Since $\nu_2(z-1)=\nu_2(z^2-1)=0$, we get that the order of $z$ modulo $\pi$ is $3$. Thus, $\nu_\pi(z^t-1)=0$ unless $3\mid t$ and $\nu_{\pi}(z^3-1)=1$. For $t=3k$, we get that 
$$
\nu_{\pi}(z^{3k}-1)=\nu_{\pi}(z^3-1)((z^3)^{k-1}+\cdots+1)=\nu_{\pi}(z^3-1)=1,
$$
provided $k$ is odd. This proves (i) when $n$ is odd. Since $z^3+1=z^3-1+2\equiv {z^3-1}\pmod {\pi^3}$, we get that $\nu_{\pi}(z^3+1)=1$. Thus, $\nu_{\pi}(z^6-1)=2$. Now $z^6+1=z^6-1+2\equiv z^6-1\pmod {\pi^3}$, so also 
$\nu_{\pi}(z^6+1)=2$. Thus, 
$$
\nu_{\pi}(z^{12}-1)=\nu_{\pi}((z^6-1)(z^6+1))=\nu_{\pi}(z^3-1)+\nu_{\pi} (z^3+1)=4.
$$
The above calculations show that $\nu_{\pi}(z^t-1)<4$ unless $12\mid t$. In particular, unless $12\mid t$, we get $\nu_2(z^t-1)=\nu_{\pi}(z^t-1)/3\le 1$, again confirming (i) in this case. However, if $12\mid t$, then setting $k:=t/12$, then since 
$$
\nu_{\pi}(z^{12}-1)=4>\frac{e_{\pi}}{p-1}(=3),
$$
we get by formula \eqref{eq:valuation} that 
$$
\nu_{\pi}(z^t-1)=\nu_{\pi}((z^{12})^k-1)=\nu_{\pi}(z^{12}-1)+\nu_{\pi}(k)=4+(\nu_{\pi}(n)-6)<\nu_{\pi}(n)=3\nu_2(n),
$$
and since $\nu_2(z^t-1)=\nu_{\pi}(z^t-1)/3$, we get the conclusion (i). 
Conclusion (ii) is even easier. Indeed, $\pi=3{\mathcal O}_{\mathbb K}$ is prime. Further, the order of $z$ in ${\mathcal O}_{\mathbb K}/\pi$ is exactly $13$ and $\nu_{3}(z^{13}-1)=1>1/12$. Thus, formula \eqref{eq:valuation} applies and gives that $\nu_3(z^t-1)=0$ unless $13\mid t$ and if $13\mid t$, then $\nu_3(z^t-1)=1+\nu_3(n)$. 
\end{proof}

\subsection{Bounds for solutions to $S-$unit equations}

It is also a well known fact from \cite{BFY} that
\begin{align}\label{2.3t}
\alpha^{n-2} < T_n < \alpha^{n-1} \quad \text{holds for all } n \geq 1.
\end{align}
Now, if $c \ge 0$ in \eqref{1.2t}, then
\begin{align*}
	2^x 3^y=T_n -c\le T_n< \alpha^{n-1}.
\end{align*}
This implies that $x\log2 +y\log 3<(n-1)\log \alpha<(n-1)\log 1.84$. Hence,
\begin{align}\label{2.4t}
	x,y<n.
\end{align}
On the other hand, if $c <0$ in \eqref{1.2t}, then
\begin{align*}
	\alpha^{n-2} < T_n =2^x 3^y+c<2^x 3^y.
\end{align*}
This implies that $\alpha^{n-2}<2^{x}3^y$, or simply,
\begin{align}\label{2.5t}
	(n-2)\log\alpha<2\max\{x\log 2, y\log 3\}.
\end{align}

The purpose of this subsection is to deduce a result from the following lemma. The following lemma is Proposition 1 from \cite{VZ}. 
\begin{lemma}[Proposition 1, \cite{VZ}]\label{lem2.8t}
	Let $\Delta>10^{80}$ be a fixed integer and assume that
	\begin{align}\label{2.10}
		2^x3^y-2^{x_1}3^{y_1}=\Delta.
	\end{align}
Then $2^x3^y<\Delta (\log\Delta)^{60\log\log\Delta}$.
\end{lemma}

We now state and prove the following consequence of Lemma \ref{lem2.8t}.
\begin{lemma}\label{lem2.9t}
	Assume that $(n,n_1,x,x_1,y,y_1)$ is a solution to $T_n-2^x3^y=T_{n_1}-2^{x_1}3^{y_1}$, with $n>310$ and $n>n_1$. Then, for $X=x\log2+y\log 3$, we have
	\[
0.08\alpha^n<\exp(X)<\alpha^n\left(n\log \alpha\right)^{60\log (n\log\alpha)},
	\]
	and
	\[
	n\log\alpha-3<X<n\log\alpha+60\left(\log (n\log\alpha)\right)^2.
	\]
\end{lemma}
\begin{proof}
	If $n>310$, then
	\[
\Delta=2^x3^y-2^{x_1}3^{y_1}=	T_n-T_{n_1}\ge T_{n-2}> T_{308}>\alpha^{306}>1.83^{306}>10^{80}.
	\]
So, we apply Lemma \ref{lem2.8t}, with $	\Delta=T_n-T_{n_1}<T_n<\alpha^{n-1}$, by \eqref{2.3t}. This yields
	\begin{align}\label{2.13t}
	\exp(X)&<\alpha^{n-1}\left(\log\alpha^{n-1}\right)^{60\log\log\alpha^{n-1}}<\alpha^{n}\left(\log\alpha^{n}\right)^{60\log\log\alpha^{n}}=\alpha^n\left(n\log \alpha\right)^{60\log (n\log\alpha)}.
	\end{align}
On the other hand, 
	\begin{align}\label{2.14t}
	0.08\alpha^n<\alpha^{n-4}<T_{n-2}\le T_n-T_{n_1}=2^x3^y-2^{x_1}3^{y_1}<2^x3^y=\exp(X).
	\end{align}
Combining \eqref{2.13t} and \eqref{2.14t}, we get
$$
0.08\alpha^n<\exp(X)<\alpha^n\left(n\log \alpha\right)^{60\log (n\log\alpha)},
$$
and taking logarithms both sides gives
\begin{align*}
n\log\alpha-3<X<n\log\alpha+60\left(\log (n\log\alpha)\right)^2,
\end{align*}
which proves Lemma \ref{lem2.9t}. 
\end{proof}

\section{Proof of Theorem \ref{1.2at}}
We go back to \eqref{1.2t} and treat the case $c=0$. This reduces the Diophantine equation to
\begin{equation}
\label{3.1t}
T_n =2^x3^y.
\end{equation} 
Now, an investigation on the largest prime factor of $k-$generalized Fibonacci numbers was made in \cite{BJJ}. Specifically, when $k=3$, it was shown in \cite{BJJ} that the only solutions to equation \eqref{3.1t} have $n\le 9$. A verification by hand in this range yields the corresponding solutions from Theorem \ref{1.2at}. \qed

\section{Proof of Theorem \ref{1.2bt}}
\subsection{An absolute upper bound on $n$}\label{subsec3.2}
For technical reasons, we assume $n>310$ and we determine an upper bound on $n$. Let  $(n,x,y)$, $(n_1,x_1,y_1)$, $(n_2, x_2,y_2)$ and $(n_3, x_3,y_3)$ be elements from $\mathbb{Z}_{\ge 0}^3$ such that
\begin{align}\label{3.2t}
	T_n-2^x3^y=T_{n_1}-2^{x_1}3^{y_1}=T_{n_2}-2^{x_2}3^{y_2}=T_{n_3}-2^{x_3}3^{y_3}.
\end{align}
Without loss of generality, we assume $n>n_1>n_2>n_3$. 
\begin{lemma}\label{lem3.1t}
	Let $c\ge 1$ be such that the Diophantine equation \eqref{1.2t} has at least two representations as 
	$$c=T_n-2^x3^y=T_{n_1}-2^{x_1}3^{y_1}.$$
	Then $$	n-n_1<9.6\cdot 10^{15} \log n.$$
\end{lemma}
\begin{proof} Since $c\ge 1$, we have $T_n>2^x3^y$ and $T_{n_1}>2^{x_1}3^{y_1}$. We go back to equation \eqref{1.2t} and rewrite it as
\begin{align*}
	T_n-2^x3^y&=T_{n_1}-2^{x_1}3^{y_1}\\
	C_\alpha\alpha^{n-1} + C_\beta\beta^{n-1}  + C_\gamma\gamma^{n-1}-2^x3^y &=C_\alpha\alpha^{n_1-1} + C_\beta\beta^{n_1-1}  + C_\gamma\gamma^{n_1-1} -2^{x_1}3^{y_1}\\
	C_\alpha\alpha^{n-1} -2^x3^y &=C_\alpha\alpha^{n_1-1} + C_\beta\beta^{n_1-1}  + C_\gamma\gamma^{n_1-1}- C_\beta\beta^{n-1}  - C_\gamma\gamma^{n-1} -2^{x_1}3^{y_1}\\
	&< C_\alpha\alpha^{n_1-1} +2<10C_\alpha\alpha^{n_1-1},
\end{align*}
for all $n>310$. So, we conclude that
\begin{align}\label{3.3t}
	\left|2^x 3^yC_\alpha^{-1}\alpha^{-(n-1)}-1\right| < 10\alpha^{-(n-n_1)}.
\end{align}
We now apply Theorem \ref{thm:Matt} on the left--hand side of \eqref{3.3t}. Let $\Gamma_0:=2^x 3^yC_\alpha^{-1}\alpha^{-(n-1)}-1=e^{\Lambda_0}-1$. Notice that $\Lambda_0\ne 0$ otherwise $\Gamma_0=0$ and we would have $C_\alpha\alpha^{n-1} =2^x3^y\in \mathbb{Z}$. Conjugating this relation by any automorphism that sends $\alpha$ to $\beta$,  we get $C_\beta\beta^{n-1} =2^x3^y$, which is a contradiction because $|C_\beta\beta^{n-1}|<1$ while $2^x3^y\ge 1$ for all $x,y\ge 0$. We use the field ${\mathbb K}:=\mathbb{Q}(\alpha)$ of
degree $D = 3$. Here, $t := 4$,
\begin{equation}\nonumber
	\begin{aligned}
		\lambda_{1}&:=2, ~~\lambda_{2}:=3, ~~\lambda_{3}:=C_\alpha, ~~\lambda_{4}:=\alpha,\\
		b_1&:=x, ~~b_2:=y, ~~b_3:=-1, ~~b_4:=-(n-1).
	\end{aligned}
\end{equation} 
Next,  $\max\{|b_1|, |b_2|, |b_3|, |b_4|\} = \max\{x,y,1,n-1\}<n$, so we can take $B:=n$. The minimal polynomial of $C_\alpha$ over the integers 
is indicated at \eqref{eq:poly}. Since its roots satisfy $|C_\alpha|$, $|C_\beta|$, $|C_\gamma|<1$, it follows that $h(C_\alpha)=\frac{1}{3}\log 44$. So, we can take
$A_1 := Dh(\lambda_{1}) = 3\log2$, $A_2 := Dh(\lambda_{2}) = 3\log3$, $A_3 := Dh(\lambda_{3}) = \log 44$, and $A_4 := Dh(\lambda_{4}) = 3\cdot\frac{1}{3} \log\alpha = \log \alpha$. 
Therefore, Theorem \ref{thm:Matt} gives,
	\begin{align}\label{3.4t}
		\log |\Gamma| &> -1.4\cdot 30^7 \cdot 4^{4.5}\cdot 3^2 (1+\log 3)(1+\log n)(3\log 2)(3\log3)(\log 44)(\log \alpha)\nonumber\\
		&> -5.6\cdot 10^{15} \log n,
	\end{align}
	where the last inequality holds for $n>310$. 
Comparing \eqref{3.3t} and \eqref{3.4t}, we get
\begin{align*}
	(n-n_1)\log \alpha-\log 10<5.6\cdot 10^{15} \log n,
\end{align*}
so that $n-n_1<9.6\cdot 10^{15}\log n$.
\end{proof}

Next, we state and prove the following result.
\begin{lemma}\label{lem3.2t}
Let	$c\ge 1$, $X := x \log 2 +y \log 3$ and $X_1:= x_1 \log 2 + y_1 \log 3$. Then
\begin{align*}
	X-X_1 <8.5\cdot 10^{30} \left(\log n\right)^2.
\end{align*}
\end{lemma}
\begin{proof}
We go back to  equation \eqref{1.2t} and rewrite it as 
\begin{align*}
	C_\alpha\alpha^{n-1} -C_\alpha\alpha^{n_1-1} -2^x3^y &= C_\beta\beta^{n_1-1}  - C_\beta\beta^{n-1}  + C_\gamma\gamma^{n_1-1}- C_\gamma\gamma^{n-1}- 2^{x_1}3^{y_1};\\ 
	\dfrac{C_\alpha\alpha^{n_1-1}\left(\alpha^{n-n_1}-1\right)}{2^x3^y}-1&=\dfrac{C_\beta\left(\beta^{n_1-1}-\beta^{n-1}\right)}{2^x3^y} + \dfrac{C_\gamma\left(\gamma^{n_1-1}-\gamma^{n-1}\right)}{2^x3^y}-\dfrac{2^{x_1}3^{y_1}}{2^x3^y},
\end{align*}
and taking absolute values, we have 
\begin{align*}
	\left|\dfrac{C_\alpha\alpha^{n_1-1}\left(\alpha^{n-n_1}-1\right)}{2^x3^y}-1\right|&\le \dfrac{1}{\exp(X)} + \dfrac{1}{\exp(X)} +  \dfrac{1}{\exp \left(X-X_1\right)}
	\le 3\exp \left(-(X-X_1)\right),
\end{align*}
where we have used the fact that the contribution of the complex conjugate zeros $\beta$ and $\gamma$ to the right--hand side of the equation above is very small, see relation \eqref{eq:approx}.

Let $\Gamma_1:=C_\alpha\alpha^{n_1-1}\left(\alpha^{n-n_1}-1\right)\cdot 2^{-x}3^{-y}-1$. Then
\begin{align}\label{3.5t}
|\Gamma_1|=\left|\dfrac{C_\alpha\alpha^{n_1-1}\left(\alpha^{n-n_1}-1\right)}{2^x3^y}-1\right|
&\le 3\exp \left(-(X-X_1)\right).	
\end{align}
Notice that $\Gamma_1=e^{\Lambda_1}-1\ne 0$, otherwise we would have
\begin{align*}
	\dfrac{C_\alpha\left(\alpha^{n-1} -\alpha^{n_1-1}\right)}{2^x3^y}=1.
\end{align*}
Taking algebraic conjugates, we would get
\begin{align*}
	1=\left|\dfrac{C_\beta\left(\beta^{n-1}  - \beta^{n_1-1}\right)}{2^x3^y}\right|<1,
\end{align*} 
a contradiction. Therefore, $\Gamma_1\ne 0$ (so, also $\Lambda_1\ne 0$).  We again use the field $\mathbb{Q}(\alpha)$ of
degree $D = 3$. Here, $t := 4$,
\begin{equation}\nonumber
	\begin{aligned}
		\lambda_{1}&:=2, ~~\lambda_{2}:=3, ~~\lambda_{3}:=\alpha, ~~\lambda_{4}:=C_{\alpha}(\alpha^{n-n_1}-1),\\
		b_1&:=-x, ~~b_2:=-y, ~~b_3:=n_1-1, ~~b_4:=1.
	\end{aligned}
\end{equation}
Next, $\max\{|b_1|, |b_2|, |b_3|,|b_4|\} = \max\{x,y,1,n_1-1\}<n$, so we can take $B:=n$. As before, we can still take
$A_1 := 3\log2$, $A_2 := 3\log3$, $A_3 := \log \alpha$. Moreover, assuming $n-n_1\ge 2$, so $\alpha^{n-n_1}-1>1$, we have 
\begin{eqnarray}
\label{eq:lambda4}
3h(\lambda_{4}) & = & 3h(C_{\alpha}(\alpha^{n-n_1}-1))\le 3h(C_{\alpha})+3h(\alpha^{n-n_1}-1)\le  \log 44+\log(\alpha^{n-n_1}-1)+2\log(1.55)\nonumber\\
& <& (n-n_1)\log \alpha+\log 44+2\log 1.55<5.7\times 10^{15}\log n.
\end{eqnarray}
where we have used Lemma \ref{lem3.1t}. 
The number $1.55$ above is an upper bound on 
$$
|\beta^{n-n_1}-1|\le 1+|\beta|^{n-n_1}|\le 1+1/\alpha<1.55
$$
and also on $|\gamma^{n-n_1}-1|$. This was if $n-n_1\ge 2$, but if $n-n_1=1$, the above bound \eqref{eq:lambda4} still holds. 
So, we take $A_4:=5.7\cdot 10^{15} \log n$. 

Then, by Theorem \ref{thm:Matt},
\begin{align}\label{3.6t}
	\log |\Gamma_1| &> -1.4\cdot 30^7 \cdot 4^{4.5}\cdot 3^2 (1+\log 3)(1+\log n)(3\log 2)(3\log3)(\log \alpha)(5.7\cdot 10^{15}\log n)\nonumber\\
	&> -8.4\cdot 10^{30} \left(\log n\right)^2.
\end{align}
Comparing \eqref{3.5t} and \eqref{3.6t}, we get
\begin{align}
	X-X_1 <8.5\cdot 10^{30} \left(\log n\right)^2.
\end{align}
This proves Lemma \ref{lem3.2t}.
\end{proof}

To proceed further, let us write $x_{min} := \min\{x,x_1,x_2\}$, $y_{min} := \min\{y,y_1,y_2\}$. We state and prove the following result.
\begin{lemma}\label{lem3.3t}
	Assume that $c\ge 1$. Then either
	\[
	x_{min}, ~y_{min} <  10^{22}(\log n)^3,
	\]
	or $n<36100$.
\end{lemma}
\begin{proof}
Again, we go back to equation \eqref{1.2t} and assume it has two solutions $(n,x,y)$, $(n_1,x_1,y_1)$. We rewrite it as	
\begin{align}\label{3.8t}
	C_\alpha\alpha^{n-1} -C_\alpha\alpha^{n_1-1}  +C_\beta\beta^{n-1}-C_\beta\beta^{n_1-1}+ C_\gamma\gamma^{n-1}-   C_\gamma\gamma^{n_1-1} &=2^x3^y- 2^{x_1}3^{y_1};\nonumber\\
	C_\alpha\alpha^{n-1}\left(1-\alpha^{n_1-n}\right)+C_\beta\beta^{n-1}\left(1-\beta^{n_1-n}\right)+C_\gamma\gamma^{n-1}\left(1-\gamma^{n_1-n}\right)&=2^x3^y- 2^{x_1}3^{y_1}.
\end{align}	
At this point, if we assume a third solution $(n_2,x_2,y_2)$ to \eqref{1.2t}, then we can also rewrite a relation analogous to \eqref{3.8t} as
\begin{align}\label{3.9t}
	C_\alpha\alpha^{n-1}\left(1-\alpha^{n_2-n}\right)+C_\beta\beta^{n-1}\left(1-\beta^{n_2-n}\right)+C_\gamma\gamma^{n-1}\left(1-\gamma^{n_2-n}\right)&=2^x3^y- 2^{x_2}3^{y_2}.
\end{align}	
We now eliminate $C_\gamma\gamma^{n-1}$ from \eqref{3.8t} and \eqref{3.9t} to get
\begin{align}\label{3.10t}
	&~~~~C_\alpha\alpha^{n-1}\left[\left(1-\alpha^{n_1-n}\right)\left(1-\gamma^{n_2-n}\right)-\left(1-\alpha^{n_2-n}\right)\left(1-\gamma^{n_1-n}\right)\right] +\nonumber\\
	&~~~~
	C_\beta\beta^{n-1}\left[\left(1-\beta^{n_1-n}\right)\left(1-\gamma^{n_2-n}\right)-\left(1-\beta^{n_2-n}\right)\left(1-\gamma^{n_1-n}\right)\right]\nonumber\\
	&=\left(2^x3^y- 2^{x_1}3^{y_1}\right)\left(1-\gamma^{n_2-n}\right)-
	\left(2^x3^y- 2^{x_2}3^{y_2}\right)\left(1-\gamma^{n_1-n}\right).
\end{align}
We first show that the left--hand side of \eqref{3.10t} is nonzero. Indeed, if it was equal to zero, then we would have $\left(2^x3^y- 2^{x_1}3^{y_1}\right)\left(1-\gamma^{n_2-n}\right)-
\left(2^x3^y- 2^{x_2}3^{y_2}\right)\left(1-\gamma^{n_1-n}\right)=0,$ so that
\begin{align*}
	\dfrac{1-\gamma^{n_1-n}}{1-\gamma^{n_2-n}} = \dfrac{2^x3^y- 2^{x_1}3^{y_1}}{2^x3^y- 2^{x_2}3^{y_2}}.
\end{align*}
Taking algebraic conjugates, we get
\begin{align}\label{3.11t}
	\dfrac{1-\gamma^{n_1-n}}{1-\gamma^{n_2-n}} =\dfrac{1-\alpha^{n_1-n}}{1-\alpha^{n_2-n}}.
\end{align}
Moreover, $n>n_1>n_2$ implies $n-n_2>n-n_1$, and Lemma \ref{lem:gamt} tells us that \eqref{3.11t} has no integer solutions $n-n_1$, $n-n_2$, with $n-n_2>n-n_1$. Hence, the left--hand side of \eqref{3.10t} is nonzero. Note also that the term $\left(1-\beta^{n_1-n}\right)\left(1-\gamma^{n_2-n}\right)-\left(1-\beta^{n_2-n}\right)\left(1-\gamma^{n_1-n}\right)$ appearing on the right--hand side of \eqref{3.10t} is nonzero. Indeed, if it was zero, then we would have 
\begin{align*}
	\dfrac{1-\gamma^{n_1-n}}{1-\gamma^{n_2-n}} =\dfrac{1-\beta^{n_1-n}}{1-\beta^{n_2-n}},
\end{align*}
a relation which is  an algebraic conjugate of \eqref{3.11t}, which cannot hold. Thus, the above relation is nonzero.

Next, we go back to \eqref{3.10t} and rewrite it as
\begin{align*}
&\dfrac{C_\alpha}{C_\beta}\left(\dfrac{\alpha}{\beta}\right)^{n-1}\dfrac{\left(1-\alpha^{n_1-n}\right)\left(1-\gamma^{n_2-n}\right)-\left(1-\alpha^{n_2-n}\right)\left(1-\gamma^{n_1-n}\right)}{\left(1-\beta^{n_1-n}\right)\left(1-\gamma^{n_2-n}\right)-\left(1-\beta^{n_2-n}\right)\left(1-\gamma^{n_1-n}\right)}+1
\\ &~~~~~~~~~~~~~~=
\dfrac{\left(2^x3^y- 2^{x_1}3^{y_1}\right)\left(1-\gamma^{n_2-n}\right)-
	\left(2^x3^y- 2^{x_2}3^{y_2}\right)\left(1-\gamma^{n_1-n}\right)}{C_\beta\beta^{n-1}\left[\left(1-\beta^{n_1-n}\right)\left(1-\gamma^{n_2-n}\right)-\left(1-\beta^{n_2-n}\right)\left(1-\gamma^{n_1-n}\right)\right]}
\\ &~~~~~~~~~~~~~~=
\dfrac{2^{x_{min}}3^{y_{min}}\cdot A}{C_\beta\beta^{n-1}\left[\left(1-\beta^{n_1-n}\right)\left(1-\gamma^{n_2-n}\right)-\left(1-\beta^{n_2-n}\right)\left(1-\gamma^{n_1-n}\right)\right]},
\end{align*}
where $$A:=\left( 2^{x-x_{min}}3^{y-y_{min}}-2^{x_1-x_{min}}3^{y_1-y_{min}} \right)\left(1-\gamma^{n_2-n}\right)-\left( 2^{x-x_{min}}3^{y-y_{min}}-2^{x_2-x_{min}}3^{y_2-y_{min}} \right)\left(1-\gamma^{n_1-n}\right).$$
Since $\nu_p(\beta)=0$, for $p=2,3$ and 
\[
\nu_p(C_\beta) = 
\begin{cases} 
	-2/3, & \text{if }~ p=2; \\
	0, & \text{if } ~p=3,
\end{cases}
\]
we have
	\begin{align*}
		&\nu_2\left(\dfrac{C_\alpha}{C_\beta}\left(\dfrac{\alpha}{\beta}\right)^{n-1}\dfrac{\left(1-\alpha^{n_1-n}\right)\left(1-\gamma^{n_2-n}\right)-\left(1-\alpha^{n_2-n}\right)\left(1-\gamma^{n_1-n}\right)}{\left(1-\beta^{n_1-n}\right)\left(1-\gamma^{n_2-n}\right)-\left(1-\beta^{n_2-n}\right)\left(1-\gamma^{n_1-n}\right)}+1 \right)\nonumber\\
		&\quad\quad=x_{min}+ \nu_2(A)-(-2/3)-\nu_2\left(\left(1-\beta^{n_1-n}\right)\left(1-\gamma^{n_2-n}\right)-\left(1-\beta^{n_2-n}\right)\left(1-\gamma^{n_1-n}\right)\right).
\end{align*}
Therefore, since $\nu_2(A)\ge 0$ because $A$ is an algebraic integer, we get
\begin{align}\label{3.12t}
 x_{min}&\le 
\nu_2\left(\dfrac{C_\alpha}{C_\beta}\left(\dfrac{\alpha}{\beta}\right)^{n-1}\dfrac{\left(1-\alpha^{n_1-n}\right)\left(1-\gamma^{n_2-n}\right)-\left(1-\alpha^{n_2-n}\right)\left(1-\gamma^{n_1-n}\right)}{\left(1-\beta^{n_1-n}\right)\left(1-\gamma^{n_2-n}\right)-\left(1-\beta^{n_2-n}\right)\left(1-\gamma^{n_1-n}\right)}+1 \right)\nonumber\\
&\quad\quad+ \nu_2\left(\left(1-\beta^{n_1-n}\right)\left(1-\gamma^{n_2-n}\right)-\left(1-\beta^{n_2-n}\right)\left(1-\gamma^{n_1-n}\right)\right).
\end{align}

In a similar way, 
\begin{align}\label{3.13t}	
y_{min}&\le 
	\nu_3\left(\dfrac{C_\alpha}{C_\beta}\left(\dfrac{\alpha}{\beta}\right)^{n-1}\dfrac{\left(1-\alpha^{n_1-n}\right)\left(1-\gamma^{n_2-n}\right)-\left(1-\alpha^{n_2-n}\right)\left(1-\gamma^{n_1-n}\right)}{\left(1-\beta^{n_1-n}\right)\left(1-\gamma^{n_2-n}\right)-\left(1-\beta^{n_2-n}\right)\left(1-\gamma^{n_1-n}\right)}+1 \right)\nonumber\\
	&\quad\quad+ \nu_3\left(\left(1-\beta^{n_1-n}\right)\left(1-\gamma^{n_2-n}\right)-\left(1-\beta^{n_2-n}\right)\left(1-\gamma^{n_1-n}\right)\right).
\end{align}	

Next, we estimate
\begin{align}\label{3.14t}
\nu_p\left(\left(1-\beta^{n_1-n}\right)\left(1-\gamma^{n_2-n}\right)-\left(1-\beta^{n_2-n}\right)\left(1-\gamma^{n_1-n}\right)\right),	
\end{align}
for $p=2,3$. Well, the number shown at \eqref{3.14t} is an algebraic integer all whose conjugates are bounded in absolute value by 
$$
2(1+\alpha^{(n-n_1)/2})(1+\alpha^{(n-n_2)/2})<\alpha^{(n-n_1)/2+(n-n_2)/2+6}\le \alpha^{n-n_2+6}.
$$
Thus, the last coefficient of the minimal polynomial for the number shown inside the valuation at \eqref{3.14t} is at most $\alpha^{6(n-n_2)+36}$. This shows that
$$
\nu_p\left(\left(1-\beta^{n_1-n}\right)\left(1-\gamma^{n_2-n}\right)-\left(1-\beta^{n_2-n}\right)\left(1-\gamma^{n_1-n}\right)\right)\le \frac{(6(n-n_2)+30)\log \alpha}{\log p}.
$$
Using Lemma \ref{3.1t}, we get that the above bounds are at most 
$$
5.7\cdot 10^{16}\log n\quad (p=2)\qquad {\text{\rm and}}\qquad 3.2\cdot 10^{16}\log n\quad (p=3).
$$
At this point, relations \eqref{3.12t} and \eqref{3.13t} become 
\begin{align}\label{3.15t}
	x_{min}&\le 
	\nu_2\left(\dfrac{C_\alpha}{C_\beta}\left(\dfrac{\alpha}{\beta}\right)^{n-1}\dfrac{\left(1-\alpha^{n_1-n}\right)\left(1-\gamma^{n_2-n}\right)-\left(1-\alpha^{n_2-n}\right)\left(1-\gamma^{n_1-n}\right)}{\left(1-\beta^{n_1-n}\right)\left(1-\gamma^{n_2-n}\right)-\left(1-\beta^{n_2-n}\right)\left(1-\gamma^{n_1-n}\right)}+1 \right)+5.7\cdot 10^{16}\log n
\end{align}
and 
\begin{align}\label{3.16t}	
	y_{min}&\le 
	\nu_3\left(\dfrac{C_\alpha}{C_\beta}\left(\dfrac{\alpha}{\beta}\right)^{n-1}\dfrac{\left(1-\alpha^{n_1-n}\right)\left(1-\gamma^{n_2-n}\right)-\left(1-\alpha^{n_2-n}\right)\left(1-\gamma^{n_1-n}\right)}{\left(1-\beta^{n_1-n}\right)\left(1-\gamma^{n_2-n}\right)-\left(1-\beta^{n_2-n}\right)\left(1-\gamma^{n_1-n}\right)}+1 \right)+ 3.2\cdot 10^{16}\log n.
\end{align}	
Lastly, we estimate
\begin{align}\label{3.17t}
\nu_p\left(\dfrac{C_\alpha}{C_\beta}\left(\dfrac{\alpha}{\beta}\right)^{n-1}\dfrac{\left(1-\alpha^{n_1-n}\right)\left(1-\gamma^{n_2-n}\right)-\left(1-\alpha^{n_2-n}\right)\left(1-\gamma^{n_1-n}\right)}{\left(1-\beta^{n_1-n}\right)\left(1-\gamma^{n_2-n}\right)-\left(1-\beta^{n_2-n}\right)\left(1-\gamma^{n_1-n}\right)}+1 \right).	
\end{align}
Let 
$$
\lambda_{1}:=\dfrac{\alpha}{\beta}\qquad {\text{\rm and}}\qquad \lambda_{2}:=-\dfrac{C_\alpha}{C_\beta}\left(\dfrac{\left(1-\alpha^{n_1-n}\right)\left(1-\gamma^{n_2-n}\right)-\left(1-\alpha^{n_2-n}\right)\left(1-\gamma^{n_1-n}\right)}{\left(1-\beta^{n_1-n}\right)\left(1-\gamma^{n_2-n}\right)-\left(1-\beta^{n_2-n}\right)\left(1-\gamma^{n_1-n}\right)}\right).
$$ 
We distinguish two cases.

\medskip

{\bf Case 1.} {\it $\lambda_1$ and $\lambda_2$ are multiplicatively independent. }

\medskip

We can apply Lemma \ref{lem:Bugt} to \eqref{3.17t} with the field $\mathbb{Q}(\alpha,\beta)$ of degree $D := 6$. Since $h(\lambda_{1})=3(\log \alpha)/3=\log\alpha$, we can choose $h(\lambda_{1}):=\log\alpha$ for all cases $p=2,3$
(note that $\log \alpha>\log p/D$ for $p=2,3$).  Further, 
\begin{align*}
	h'(\lambda_{2}) &\le 2h(C_\alpha)+ 4h\left(1-\alpha^{n_1-n}\right)+4\log(1-\alpha^{n-n_2})+2\log 2\\
	&< \frac{2\log 44}{3}+2\log 2+\frac{4}{3}\left(\log(1+|\beta|^{n_1-n})+\log(1+|\gamma|^{n_2-n})\right)\\
	& <  4+\frac{4}{3}\left(\log(\alpha^{(n-n_1)/2}+1)+\log(\alpha^{(n-n_2)/2}+1)\right)\\
	& <  4+\frac{4}{3}\left((n-n_1)/2+2+(n-n_2)/2+2)\right)\log \alpha\\
	& <  4+\frac{8(n-n_2+4)}{3}\log \alpha<1.6\times 10^{16}\log n,
\end{align*}
where we have used Lemma \ref{lem3.1t} and $n>310$. Therefore,
\[
E'=\dfrac{b_1}{h'(\lambda_{2})}+\dfrac{b_2}{h'(\lambda_{1})}=\dfrac{n-1}{h'(\lambda_{1})}+\dfrac{1}{h'(\lambda_{2})}<n.
\]
We now have
\begin{align*}
	E&=\max\left\{\log E'+\log\log p+0.4, 10, 10\log p\right\}\\
	&=\max\left\{\log n +\log\log p+0.4, 10,10\log p\right\}.
\end{align*}
If $n>36100$, then $E<1+\log n$ in both cases $p=2,3$. Therefore, Lemma \ref{lem:Bugt} gives
\begin{align*}
&\nu_p\left(\dfrac{C_\alpha}{C_\beta}\left(\dfrac{\alpha}{\beta}\right)^{n-1}\dfrac{\left(1-\alpha^{n_1-n}\right)\left(1-\gamma^{n_2-n}\right)-\left(1-\alpha^{n_2-n}\right)\left(1-\gamma^{n_1-n}\right)}{\left(1-\beta^{n_1-n}\right)\left(1-\gamma^{n_2-n}\right)-\left(1-\beta^{n_2-n}\right)\left(1-\gamma^{n_1-n}\right)}+1 \right)\\
&\le \dfrac{24pg}{(p-1)(\log p)^4} E^2D^4   h'(\lambda_{1})h'(\lambda_{2}  )\\
	&<\dfrac{24pg}{(p-1)(\log p)^4} (\log n)^2\left(1+\frac{1}{\log 36100}\right)^2\cdot 6^4  ( \log\alpha )\cdot 1.6\cdot10^{16}\log n\\
	&<\dfrac{3.7\cdot 10^{20}pg}{(p-1)(\log p)^4} \left(\log n\right)^3.
\end{align*}
Hence, inequalities \eqref{3.15t} and \eqref{3.16t} become 
\begin{align*}
	x_{min}&< 
	\dfrac{3.7\cdot 10^{20}\cdot 2\cdot 3}{(2-1)(\log 2)^4} \left(\log n\right)^3+ 5.7\times 10^{16} \log n<10^{22} (\log n)^3
\end{align*}
and 
\begin{align*}
	y_{min}&<\dfrac{3.7\cdot 10^{20}\cdot 3\cdot 13}{(3-1)(\log 3)^4} \left(\log n\right)^3+ 3.2\times 10^{16} \log n<5\times 10^{21} (\log n)^3.
\end{align*}	
In the above, we used $g=3$ when $p=2$ and $g=13$ when $p=3$ as in the proof of Lemma \ref{zLucas}. 
This was in case $\lambda_1$ and $\lambda_2$ are multiplicatively independent. 

\medskip

{\bf Case 2.} {\it $\lambda_1$ and $\lambda_2$ are multiplicatively dependent.}

\medskip

In this case, $\lambda_1=z^3$ and $\lambda_2=\pm z^{t}$ for some integer $t$ with $|t|\le 80n+244$ by Lemma \ref{lem:a/b}. Thus, 
$$
\lambda_1^{n-1}\lambda_2-1=\pm z^{t+3(n-1)}+1\mid z^{2t+6n-2}-1.
$$
The exponent of $z$ above is in absolute value is at most $2|t|+6n-2\le 186n+500<n^2$. Lemma \ref{zLucas} now shows that 
$$
\nu_p(\lambda_1^{n-1}\lambda_2+1)\le 1+\frac{\log n^2}{\log p}<3\log n
$$
for $n>36100$ and $p\in \{2,3\}$. Comparing this to Case 1, we see that the bounds from Case 1 still hold. 
This completes the proof of Lemma \ref{lem3.3t}.
\end{proof}

Lastly, we consider a fourth solution $(n_3,x_3,y_3)$, with $n>n_1>n_2>n_3$ and we find an absolute bound for $n$. We prove the following result.
\begin{lemma}\label{lem3.4t}
	If $c\ge 1$ and $n>36100$,  then
	\[
	n<1.2\cdot 10^{37},\qquad x<1.1\cdot 10^{37}, \qquad y<6.8\cdot 10^{36} .
	\]
\end{lemma}
\begin{proof}
Lemma \ref{lem3.3t} indicates that within any set of three solutions, the smallest values of $x$ and $y$ are constrained to be less than $10^{22}(\log n)^3$. Consequently, in a quadrat of solutions, it is possible for at most one solution to have an $x$--value exceeding $10^{22}(\log n)^3$
and similarly, at most one solution can have a $y$--value surpassing this bound. Therefore, it follows that at least one solution in each set will have both $x$ and $y$ confined within this limit. This particularly implies that the smallest solution adheres to these bounds. Hence,
	$$
	X_3=x_3\log 2+y_3\log 3<(\log 2+\log 3)\cdot  10^{22}(\log n)^3<1.8\cdot 10^{22} (\log n)^3.
	$$
	With Lemmas \ref{lem2.9t} and \ref{lem3.2t}, we get
	\begin{align*}
		n\log\alpha-3&<X\\
		&=X_3+(X_2-X_3)+(X_1-X_2)+(X-X_1)\\
		&< 1.8\cdot10^{22}(\log n)^3 +3\cdot 8.5\cdot 10^{30} \left(\log n\right)^2\\
		&<2.44\cdot 10^{30} \left(\log n\right)^3,
	\end{align*} 	
	which implies
	\begin{align}\label{3.18t}
		\dfrac{n}{(\log n)^3}<4.1\cdot 10^{30}.
	\end{align}	
	We apply Lemma \ref{Guzt} to inequality \eqref{3.18t} above with $ z:=n $, $ m:=3 $, $T:=4.1\cdot10^{30}$.
	Since $T>(4\cdot 3^2)^3=46656$, we get
	$$n<2^m T(\log T)^m = 2^3 \cdot 4.1\cdot 10^{30}(\log 4.1\cdot 10^{30})^3 < 1.2\cdot 10^{37}.$$	
	Further, we have by Lemma \ref{lem2.9t} that
	\begin{align*}
		X&<n\log\alpha+60\left(\log (n\log\alpha)\right)^2;\\
		x\log2+y\log3&<1.2\cdot 10^{37}\log\alpha+60\left(\log (1.2\cdot 10^{37}\log\alpha)\right)^2<7.4\cdot 10^{36}.
	\end{align*}	
	This gives $x<1.1\cdot 10^{37}$ and $y<6.8\cdot 10^{36}$ and completes the proof of Lemma \ref{lem3.4t}.
\end{proof}

\subsection{Reduction of the upper bound on $n$} \label{subsec3.2.2}
Here, we use the LLL--reduction method, the theory of continued fractions and as well as $p-$adic reduction methods due to \cite{PET} to find a rather small bound for $n$.

To begin, we go back to equation \eqref{3.3t}. Assuming $n-n_1\ge 5$, we can write
\begin{align*}
	|\Lambda_0|=\left|\log C_\alpha-x\log 2-y\log 3+(n-1)\log\alpha\right|\le 15\alpha^{-(n-n_1)},
\end{align*}
where we used Lemma \ref{lem2.1t} with $n-n_1\ge 5$ since $\alpha^{n-n_1}\ge \alpha^5>2$. So, we consider the approximation lattice
$$ \mathcal{A}=\begin{pmatrix}
	1 & 0  & 0 \\
	0 & 1 & 0 \\
	\lfloor M\log 2\rfloor & \lfloor M\log 3\rfloor& \lfloor M\log\alpha \rfloor
\end{pmatrix},$$
with $M:= 10^{112}$ and choose $v:=\left(0,0,-\lfloor M\log C_\alpha \rfloor\right)$. Now, by Lemma \ref{lem2.5t}, we get $$|\Lambda|\ge c_1:=10^{-40} \qquad\text{and}\qquad c_2:=2.15\cdot 10^{38}.$$
Moreover, by Lemma \ref{lem3.4t}, we have
\[
x<A_1:=1.1\cdot 10^{37},~~y<A_2:=6.8\cdot 10^{36},~~n-1<A_3:=1.2\cdot 10^{37}.
\]
So, Lemma \ref{lem2.6t} gives $S=1.7\cdot 10^{74}$ and $T=1.49\cdot 10^{37}$. Since $c_2^2\ge T^2+S$, then choosing $c_3:=15$ and $c_4:=\log\alpha$, we get $n-n_1\le 285$.

Next, we now go back to equation \eqref{3.5t}. Assume that $X-X_1\ge 2$. We can then write
\begin{align*}
	|\Lambda_1|=\left|(n_1-1)\log\alpha+\log C_\alpha\left(\alpha^{n-n_1}-1\right)-x\log 2-y\log 3 \right|
	&< 4.5\exp \left(-(X-X_1)\right),
\end{align*}
where we used Lemma \ref{lem2.1t} together with the fact that $\exp(X-X_1)\ge \exp(2)>6$. So, we use the same approximation lattice
$$ \mathcal{A}=\begin{pmatrix}
	1 & 0  & 0 \\
	0 & 1 & 0 \\
	\lfloor M\log 2\rfloor & \lfloor M\log 3\rfloor& \lfloor M\log\alpha \rfloor
\end{pmatrix} ,$$
with $M:= 10^{113}$ and choose $v:=\left(0,0,-\lfloor M\log C_\alpha\left(\alpha^{n-n_1}-1\right) \rfloor\right)$. It turns out that for all values $1\le n-n_1\le 285$, the chosen constant $M$ is sufficiently large, so we can still apply Lemma \ref{lem2.6t}. By Lemma \ref{lem2.5t}, we maintain the lower bound 
$|\Lambda_1|\ge c_1:=10^{-40}$, $c_2:=6.13\cdot10^{37}$
and by Lemma \ref{lem3.4t}, we also have 
\[
x<A_1:=1.1\cdot 10^{37},~~y<A_2:=6.8\cdot 10^{36},~~n-1<A_3:=1.2\cdot 10^{37}.
\]
So, Lemma \ref{lem2.6t} still gives the same values of $S$ and $T$ as before. Since $c_2^2\ge T^2+S$, we now choose $c_3:=4.5$ and $c_4:=1$ and we get $X-X_1\le 174$.

Next, we find reduced bounds on $x_{min}$ and $y_{min}$ using $p-$adic reduction methods due to \cite{PET}. We go back to equation \eqref{1.2t} and assume it has two solutions $(n,x,y)$, $(n_1,x_1,y_1)$. We rewrite it as 
$$c=T_n-2^x3^y=T_{n_1}-2^{x_1}3^{y_1},$$
so that we can determine $\nu_p\left(T_n-T_{n_1}\right)$. Notice that from the above equation, 
\begin{align*}
	\nu_p\left(T_n-T_{n_1}\right)&=\nu_p\left(2^x3^y-2^{x_1}3^{y_1}\right)\\
	&=\nu_p\left(2^{x_{min}}3^{y_{min}}\left(2^{x-x_{min}}3^{y-y_{min}}-2^{x_1-x_{min}}3^{y_1-y_{min}}\right)\right)\\
	&=\nu_p\left(2^{x_{min}}3^{y_{min}}\right)+\nu_p\left(2^{x-x_{min}}3^{y-y_{min}}-2^{x_1-x_{min}}3^{y_1-y_{min}}\right),
\end{align*}
and thus
\begin{align}\label{3.19t}
\nu_p\left(2^{x_{min}}3^{y_{min}}\right)
	&=	\nu_p\left(T_n-T_{n_1}\right)-\nu_p\left(2^{x-x_{min}}3^{y-y_{min}}-2^{x_1-x_{min}}3^{y_1-y_{min}}\right)\nonumber\\
	&\le \nu_p\left(T_n-T_{n_1}\right) .
\end{align}
Now, we determine $\nu_p\left(T_n-T_{n_1}\right) $ for all $n-n_1\ge 1$, $n<1.2\cdot 10^{37}$ and $p=2,3$. Here, we show how to do it for $p=2$, and then we automate the process in SageMath. Since $n<1.2\cdot 10^{37}<2^{124}$, then $n$ has at most 124 binary digits. Let $d:=n-n_1\le 285$, by the results in the reduction above. So, we need an upper bound for 
$$ \nu_2 \left(T_n-T_{n_1}\right)= \nu_2 \left(T_{n_1+d}-T_{n_1}\right),$$
but since $n_1<n$, we instead bound 
$$ \nu_2 \left(T_{n+d}-T_{n}\right),\qquad d\in [1,285],\qquad n<1.2\cdot 10^{37}.$$
The Tribonacci sequence is periodic modulo $2^k$ with period $2^{k+1}$. In particular, $T_{n+d}-T_{n}$ is periodic modulo $2^8$ with period $2^9<1000$. With a simple program in SageMath, we checked over all $d\in [1,285]$ for which there is $n\le 1000$ such that $2^{8}\mid T_{n+d}-T_n$. There are $214$ such $d$'s, namely
$\{1, 3, 4, 5, 7, 8, 9, 11, 12, 13, 15, 16, 17, 19, 20, 21, 23, 24, 25, 27, 28,,\ldots,281,283,284,285\}$.
Some numbers are missing from this list. This means that for the missing numbers $d$, we have $\nu_2(T_{n+d}-T_n)\le 7$ always. 

Here, we will work out one $d$ only, for explanation. Namely, we take $d=9$. We calculate $n_0(d)\in [1, 2^9]$ such that for $n=n_0(d)$ we have that $\nu_2(L_{n+d}-L_n)\ge 8$. This is unique in this case and it is $n_0(d)=167$. So, from now on, every $n< 1.2\cdot 10^{37}$ such that $\nu_2(T_{n+d}-T_n)\ge 8$ is of the form $n=167+ 2^9 z$, for some $z\in \mathbb{Z}$. So now, we need to find out $z$ such that $\nu_2(T_{n+9}-T_n)$ is as large as possible. For this, we go to the Binet formula and get
\begin{align*}
	T_{n+9}-T_n&=\left( C_\alpha\alpha^{n+8} + C_\beta\beta^{n+8}  + C_\gamma\gamma^{n+8}\right)-\left( C_\alpha\alpha^{n-1} + C_\beta\beta^{n-1}  + C_\gamma\gamma^{n-1}\right)\\
	&=\left(\alpha^{9}-1\right)C_\alpha\alpha^{n-1}+
	\left(\beta^{9}-1\right)C_\beta\beta^{n-1}+
	\left(\gamma^{9}-1\right)C_\gamma\gamma^{n-1}\\
	&=\left(\alpha^{9}-1\right)C_\alpha\alpha^{167+ 2^9 z-1}+
	\left(\beta^{9}-1\right)C_\beta\beta^{167+ 2^9 z-1}+
	\left(\gamma^{9}-1\right)C_\gamma\gamma^{167+ 2^9 z-1}\\
	&=\left(\alpha^{9}-1\right)C_\alpha\alpha^{166}\alpha^{2^9 z}+
	\left(\beta^{9}-1\right)C_\beta\beta^{166}\beta^{2^9 z}+
	\left(\gamma^{9}-1\right)C_\gamma\gamma^{166}\gamma^{2^9 z}\\
	&=\left(\alpha^{9}-1\right)C_\alpha\alpha^{166}\exp_2\left(2^9z\log_2\alpha\right)+
	\left(\beta^{9}-1\right)C_\beta\beta^{166}\exp_2\left(2^9z\log_2\beta\right)\\&~+
	\left(\gamma^{9}-1\right)C_\gamma\gamma^{166}\exp_2\left(2^9z\log_2\gamma\right).
\end{align*}
We can apply a simple cosmetic to the above relation so that we have terms in $\log \alpha^4$, $\log \beta^4$ and $\log \gamma^4$ inside the respective exponentials. We write
\begin{align*}
	T_{n+9}-T_n
	&=\left(\alpha^{9}-1\right)C_\alpha\alpha^{166}\exp_2\left(2^7z\log_2\alpha^4\right)+
	\left(\beta^{9}-1\right)C_\beta\beta^{166}\exp_2\left(2^7z\log_2\beta^4\right)\\&~+
	\left(\gamma^{9}-1\right)C_\gamma\gamma^{166}\exp_2\left(2^7z\log_2\gamma^4\right).
\end{align*}
Moreover, 
\begin{equation}
	\label{eq:1t}
	\log_2 \alpha^4=\log_2 \left(1-\left(1-\alpha^4\right)\right)=-\sum_{n\ge 1} \frac{\left(1-\alpha^4\right)^n}{n},
\end{equation}
and in the right--hand side, 
$
|(1-\alpha^4)^{n}/n|_2\le 2^{-(n-\log n/\log 2)}<2^{-4n/3+\log n/\log 2},
$
where we have used the fact that $\nu_2(\alpha^4-1)=\nu_2((\alpha-1)(\alpha+1)(\alpha^2+1))=4/3$. This therefore shows that series given in the right--hand side of \eqref{eq:1t} converges.
For the argument in the exponential, we have 
$$
\nu_2\left(2^7 z\log_2 \alpha^4\right)\ge \nu_2(2^8z)\ge 8,
$$
so $\left|2^7z\log_2 \alpha^4\right|_2\le 2^{-8}<2^{-1}$, therefore the exponential in this input is convergent $2$--adically. The same arguments work with $\alpha$ replaced by $\beta$ or $\gamma$. We now stop the argument of the logarithm at $n=120$, so put
\begin{equation}
	\label{eq:Pt}
	P:=-\sum_{n=1}^{120} \frac{\left(1-\alpha^4\right)^n}{n},
\end{equation}
such that
$$
\log_2 \alpha^4=P-\sum_{n\ge 121} \frac{\left(1-\alpha^4\right)^n}{n}.
$$
Since
$n-\log_2 n\ge n-\log n/\log 2$ and the function $n-\log n/\log 2$ is at least $121$ for all $n\ge 128$, then $n-\nu_2(n)\ge 121$ for all $n\ge 121$. Thus, 
$\log_2 \alpha^4=P+u$,
where $\nu_2(u)\ge 121$. We therefore have,
$$
2^7 z \log_2 \alpha^4=2^7 z P+2^7 z u,
$$
so that
$$
\exp_2(2^7 z\log_2 \alpha^4)=\exp_2(2^7 zP+2^7 zu)=\exp_2(2^7zP)\exp_2(2^7 zu).
$$
We have
$$
\exp_2(y)=1+y+\frac{y^2}{2}+\cdots+\frac{y^n}{n!}+\cdots.
$$
For $\nu_2(y)\ge 2$ and $n\ge 2$ we have 
$$
\nu_2\left(\frac{y^n}{n!}\right)=n\nu_2(y)-\nu_2(n!)\ge n\nu_2(y)-(n-\sigma_2(n))>n(\nu_2(y)-1)\ge \nu_2(y),
$$
where the last inequality holds as it is equivalent to $\nu_2(y)\ge n/(n-1)$, which is so since $\nu_2(y)\ge 2\ge n/(n-1)$ for all $n\ge 2$. In the above, $\sigma_2(n)$ is the sum of the digits of $n$ in base $2$. It then follows that 
$
\exp_2(y)\equiv 1\pmod {2^{\nu_2(y)}},
$
provided $\nu_2(y)\ge 2$. Hence,
$
\exp_2(2^7 zu)\equiv 1\pmod {2^{7+\nu_2(u)}}\equiv 1\pmod {2^{128}}.
$
Thus
$$
\exp_2(2^7 z\log_2 \alpha)\equiv \exp_2(2^7 z P)\pmod {2^{128}}\equiv \sum_{k\ge 0} \frac{(2^7 z P)^k}{k!}\pmod {2^{128}}.
$$
Indeed,
$$
\nu_2\left(\frac{(2^7 z P)^k}{k!}\right)=k\nu_2(2^7 z P)-\nu_2(k!)\ge (7+\nu_2(P))k-(k-\sigma_2(k))> 7k,
$$
since $\sigma_2(k)\ge 1$ and $\nu_2(P)\ge 1$, so it follows that the above numbers are $\ge 7\cdot 19=133>128$ for $k\ge 19$. Thus, we may truncate the series at $k=18$ and write
$$
\exp_2\left(2^7 z\log_2\alpha^4\right)=\sum_{k=0}^{18} \frac{(2^7 zP)^k}{k!}\pmod {2^{128}}.
$$
The same argument works with $\alpha$ replaced by $\beta$ or $\gamma$, so we may write 
\begin{equation}
	\label{eq:QR}
	Q:=-\sum_{n=1}^{120} \frac{\left(1-\beta^4\right)^n}{n}\qquad\text{and}\qquad R:=-\sum_{n=1}^{120} \frac{\left(1-\gamma^4\right)^n}{n},
\end{equation}
so that
$$
\exp_2(2^7 z \log_2 \beta^4)=\sum_{k=0}^{18} \frac{(2^7 zQ)^k}{k!}\pmod {2^{128}}\quad\text{and}\quad
\exp_2(2^7 z \log_2 \gamma^4)=\sum_{k=0}^{18} \frac{(2^7 zR)^k}{k!}\pmod {2^{128}}.
$$
Thus,
\begin{align*}
T_{n+9}-T_n=\sum_{k=0}^{18} \frac{\left(\alpha^{9}-1\right)C_\alpha\alpha^{166} (2^7 zP)^k+\left(\beta^{9}-1\right)C_\beta\beta^{166}(2^7 z Q)^k+\left(\gamma^{9}-1\right)C_\gamma\gamma^{166}(2^7 z R)^k}{k!}\pmod {2^{128}}.
\end{align*}
The right--hand side above is a polynomial of degree $18$ in $z$ whose coefficients are rational numbers which are $2$--adic integers (that is, the numerators of those rational numbers are always odd). We will show that in our range the above expression is never $0$ modulo $2^{128}$. This will show that
$\nu_2(T_{n+9}-T_n)<128$ for $n<1.2\cdot 10^{37}.$

We need to find these numbers which is not so easy in SageMath as $P$, $Q$ and $R$ involve large powers of $\alpha$, $\beta$ and $\gamma$ respectively. Nevertheless, we can compute $A:=P+Q+R$, $B:=PQ+PR+QR$ and  $C:=PQR$. Next, the coefficients 
\begin{align}\label{3.23t}
u_k:=\left(\alpha^{9}-1\right)C_\alpha\alpha^{166} P^k+\left(\beta^{9}-1\right)C_\beta\beta^{166}Q^k+\left(\gamma^{9}-1\right)C_\gamma\gamma^{166}R^k
\end{align}
form a linearly recurrence sequence of recurrence 
$$
u_{k+3}=Au_{k+2}-Bu_{k+1}+Cu_k,\qquad {\text{\rm for}}\qquad k\ge 0,
$$
with $u_0$, $u_1$ and $u_2$ obtained from \eqref{3.23t} when $k=0,1,2$ respectively. So, we can compute all the remaining ones iteratively and look at the polynomial
$$
f(z):=\sum_{k=0}^{18} (2^7 z)^k \frac{u_k}{k!}\pmod {2^{128}}.
$$
All coefficients \( u_k/k! \) are 2--adic integers, and we can reduce them modulo \( 2^{128} \). As a result, we obtain a polynomial in \( \mathbb{Z}/(2^{128}\mathbb{Z})[z] \). Our goal is to find a \( z \) such that this polynomial \(= 0 \pmod{2^{128}} \). This is achieved incrementally. Specifically, by reducing \( f(z) \) modulo \( 2^{10}, 2^{11}, 2^{12}, 2^{13}, \) and so forth, we deduce the necessary digit of \( z \) in the subsequent power of 2 (either 0 or 1) to ensure divisibility by higher and higher powers of 2, in accordance with Hensel's lemma. This process yields 
\[
z = 2^0 + 2^1 + 2^2 + 2^3 + 2^4 + 2^5 + 2^6 + 2^7 + 2^8 + 2^9 + \cdots,
\]
and upon extracting digits up to \( 2^{101} \), we express 
\[
z = 2^0 + 2^1 + 2^2 + 2^3 + 2^4 + 2^5 + 2^6 + 2^7 + 2^8 + 2^9 + \cdots + 2^{100} + 2^{101} + 2^{102} t,
\]
and reduce \( f(z) \) modulo \( 2^{127} \), obtaining 
\[
2^{123} (8 + 9 t) \pmod{2^{127}}.
\]
Notice that we should choose \( z \) as a multiple of 8, leading to 
\[
n \geq 2^{19} \cdot (\cdots + 2^{102} \cdot 8) = 2^{124} > 1.2 \cdot 10^{37}.
\]
This analysis demonstrates that effectively, \( \nu_2(T_{n+9} - T_n) < 126 \).

A similar analysis was performed for the other values of $d$ as well as for $p=3$. In the case of $p=3$, the period of $(T_n)_{n\ge 0}$ modulo $3^{k+1}$ is $13\cdot 3^k$, so we work $3-$adically with $\log_3 (\alpha^{13})$, $\log_3 (\beta^{13})$ and $\log_3 (\gamma^{13})$. In all cases, we obtained that $\nu_p(T_{n}-T_{n_1})<128$. Hence, in all cases, we conclude from \eqref{3.19t} that
\[x_{min},~ y_{min}< 128 .\]
	
Lastly, we find a smaller upper bound for $n$. If we write $c_X$ for the upper bound of $X-X_1$, then 
\begin{align*}
	X=X_2+(X_1-X_2)+(X-X_1)&<x_{min}\log 2 +y_{min}\log 3+2c_X,\\
	x\log 2 +y\log 3&<128\log 2 +128\log 3+2\cdot 174<578.
\end{align*}
Hence, $x<834$ and $y<527$. On the other hand, Lemma \ref{lem2.9t} implies that $n\log\alpha-3<X<578,$
so that $n<962$.

\subsection{Conclusion}
To finalize with the proof of Theorem \ref{1.2bt}, we note that for \( n > 310 \), the bounds are \( n \leq 962 \), \( x < 834 \), and \( y < 527 \). To efficiently handle large \( T_n \) values, our SageMath 9.5 code utilized batch processing, iterating through all \( (n, x, y) \) combinations within these ranges. This approach selected \( c\ge 1 \) values with at least four representations of the form \( T_n - 2^x 3^y \), aligning with the solutions in Theorem \ref{1.2bt}. The computation, performed on an 8GB RAM laptop, was completed in about 3 hours. \qed

\section{Proof of Theorem \ref{1.2ct}}
\subsection{An absolute upper bound on $n$}
We proceed as in Subsection \ref{subsec3.2} and determine an absolute upper bound on $n$. Without loss of generality, we may assume $n>n_1>310$, and continue by proving a series of results.

\begin{lemma}\label{lem3.5t}
	Let $c\le -1$ such that $(n,x,y)$ and $(n_1,x_1,y_1)$ satisfy \eqref{1.2t} with $n>n_1>310$, then
	\[	X-X_1<5.55\cdot 10^{15} \log n. \]
\end{lemma}
\begin{proof}
	Since $c\le -1$, then $T_n-2^x3^y\le -1$ and hence $2^x3^y-T_n>0$. We go back to \eqref{1.2t} and rewrite it as
	\begin{align}\label{3.24t}
		0<	2^x3^y-T_n&=2^{x_1}3^{y_1}-T_{n_1};\nonumber\\
		2^x3^y-C_\alpha\alpha^{n-1}&=2^{x_1}3^{y_1}-C_\alpha\alpha^{n_1-1}+C_\beta\left(\beta^{n-1}-\beta^{n_1-1}\right)+C_\gamma\left(\gamma^{n-1}-\gamma^{n_1-1}\right)\nonumber\\
		&<2^{x_1}3^{y_1}+C_\beta\left(\beta^{n-1}-\beta^{n_1-1}\right)+C_\gamma\left(\gamma^{n-1}-\gamma^{n_1-1}\right);\nonumber\\
		|2^x3^y-C_\alpha\alpha^{n-1}|&<2^{x_1}3^{y_1}+2,
	\end{align}
where the last inequality holds for all $n>n_1>310$. So, dividing through both sides of \eqref{3.24t} by $2^x3^y$, we conclude that
	\begin{align}\label{3.25t}
		\left|2^{-x} 3^{-y}C_\alpha\alpha^{n-1}-1\right| \le 3\exp\{-(X-X_1)\}.
	\end{align}
We now apply Theorem \ref{thm:Matt} on the left--hand side of \eqref{3.25t}. Let $\Gamma_2:=2^{-x} 3^{-y}C_\alpha\alpha^{n-1}-1=e^{\Lambda_2}-1$. Like before, $\Lambda_2\ne 0$ otherwise we would have $C_\alpha\alpha^{n-1} =2^x3^y\in \mathbb{Z}$. Conjugating this relation by any automorphism that swaps $\alpha$ and $\beta$, we get $C_\beta\beta^{n-1} =2^x3^y$, which is a contradiction because $|C_\beta\beta^{n-1}|<1$ while $2^x3^y\ge 1$ for all $x,y\ge 0$. We use the field ${\mathbb K}:=\mathbb{Q}(\alpha)$ of
degree $D = 3$. Here, $t := 4$,
\begin{equation}\nonumber
	\begin{aligned}
		\lambda_{1}&:=2, ~~\lambda_{2}:=3, ~~\lambda_{3}:=C_\alpha, ~~\lambda_{4}:=\alpha,\\
		b_1&:=-x, ~~b_2:=-y, ~~b_3:=1, ~~b_4:=n-1.
	\end{aligned}
\end{equation} 
We have $\max\{|b_1|, |b_2|, |b_3|, |b_4|\} = \max\{x,y,1,n-1\}=n$. So, we can take $B:=n$. Again, $h(C_\alpha)=\frac{1}{3}\log 44$. So, we can take
$A_1 := 3\log2$, $A_2 := 3\log3$, $A_3 := \log 44$, and $A_4 : = \log \alpha$. 
Therefore, Theorem \ref{thm:Matt} gives,
\begin{align}\label{3.26t}
	\log |\Gamma_2| &> -1.4\cdot 30^7 \cdot 4^{4.5}\cdot 3^2 (1+\log 3)(1+\log n)(3\log 2)(3\log3)(\log 44)(\log \alpha)\nonumber\\
	&> -5.5\cdot 10^{15} \log n,
\end{align}
where the last inequality holds for $n>310$. 
Comparing \eqref{3.25t} and \eqref{3.26t}, we get
	\begin{align}\label{5.4h}
		X-X_1<5.55\cdot 10^{15} \log n.
	\end{align}
	This proves Lemma \ref{lem3.5t}.
\end{proof}

Next, we prove the following.
\begin{lemma}\label{lem3.6t}
	Let $c\le -1$ such that $(n,x,y)$, $(n_1,x_1,y_1)$ and  $(n_2,x_2,y_2)$ satisfy \eqref{1.2t} with $n>n_1>n_2$ and $n_1>310$, then
	\[	n-n_1<2\cdot 10^{15}(\log n)^2. \]
\end{lemma}
\begin{proof}
	Since we assume a third solution $(n_2,x_2,y_2)$ to \eqref{1.2t}, then Lemma \ref{lem2.9t} also holds for $n_1>n_2$. By hypothesis, $n_1>310$, so 
	\[
	0.08\alpha^{n_1}<\exp(X_1)=2^{x_1}3^{y_1}<\alpha^{n_1}\left(n_1\log \alpha\right)^{60\log (n_1\log\alpha)}.
	\]
	We can then rewrite \eqref{3.24t} with the above inequality and get
\begin{align*}
\left|2^x3^y-C_\alpha\alpha^{n-1}\right|&<2^{x_1}3^{y_1}+2\\
		&<\alpha^{n_1}\left(n_1\log \alpha\right)^{60\log (n_1\log\alpha)}+2\\
		&<1.01\alpha^{n_1}\left(n_1\log \alpha\right)^{60\log (n_1\log\alpha)}.
	\end{align*}
	Dividing through by $C_\alpha\alpha^{n-1}$, we obtain
	\begin{align}\label{5.5h}
		\left|2^{x} 3^{y}C_\alpha^{-1}\alpha^{-(n-1)}-1\right| &\le \dfrac{1.01\alpha}{C_\alpha}\alpha^{-(n-n_1)}\left(n_1\log \alpha\right)^{60\log (n_1\log\alpha)}\nonumber\\
		&<3\alpha^{-(n-n_1)}\left(n_1\log \alpha\right)^{60\log (n_1\log\alpha)}.
	\end{align}
Let $\Gamma_3:=2^{x} 3^{y}C_\alpha^{-1}\alpha^{-(n-1)}-1=e^{\Lambda_3}-1$. By the same arguments and data used in the proof of Lemma \ref{lem3.5t} above, we conclude by Matveev's Theorem \ref{thm:Matt} that
	\begin{align*}
		(n-n_1)\log\alpha-\log 3 -60\left[\log (n_1\log\alpha)\right]^2&<5.5\cdot 10^{15} \log n;\\
		n-n_1&<\dfrac{60\left[\log (n_1\log\alpha)\right]^2+\log 3 +5.5\cdot 10^{15} \log n}{\log\alpha}\\
		&<2\cdot 10^{15}(\log n)^2,
	\end{align*}
where we have used the fact that $n>n_1>310$. Hence, Lemma \ref{lem3.6t} is proved.
\end{proof}

Next, we retain the notation $x_{min}:=\min\{x,x_1,x_2\}$ and $y_{min}:=\min\{y,y_1,y_2\}$ and prove the following result.
\begin{lemma}\label{lem3.7t}
	Let $c\le -1$ such that $(n,x,y)$, $(n_1,x_1,y_1)$, $(n_2,x_2,y_2)$ and  $(n_3,x_3,y_3)$ satisfy \eqref{1.2t} with $n>n_1>n_2>n_3$ and $n_1>310$, then
	either	\[	x_{min},~y_{\min}<  2\cdot 10^{21}\left(\log n\right)^4, \]
	or $n<36100$.
\end{lemma}
\begin{proof}
We follow the same arguments as in the proof of Lemma \ref{lem3.3t} but using the bounds from Lemma \ref{lem3.6t}. In particular, we consider again the $p-$adic valuations \eqref{3.12t} and \eqref{3.13t}.
We have,
\begin{align}\label{3.29t}
	x_{min}&\le 
	\nu_2\left(\dfrac{C_\alpha}{C_\beta}\left(\dfrac{\alpha}{\beta}\right)^{n-1}\dfrac{\left(1-\alpha^{n_1-n}\right)\left(1-\gamma^{n_2-n}\right)-\left(1-\alpha^{n_2-n}\right)\left(1-\gamma^{n_1-n}\right)}{\left(1-\beta^{n_1-n}\right)\left(1-\gamma^{n_2-n}\right)-\left(1-\beta^{n_2-n}\right)\left(1-\gamma^{n_1-n}\right)}+1 \right)\nonumber\\
	&\quad\quad+ \nu_2\left(\left(1-\beta^{n_1-n}\right)\left(1-\gamma^{n_2-n}\right)-\left(1-\beta^{n_2-n}\right)\left(1-\gamma^{n_1-n}\right)\right),
\end{align}
and
\begin{align}\label{3.30t}	
	y_{min}&\le 
	\nu_3\left(\dfrac{C_\alpha}{C_\beta}\left(\dfrac{\alpha}{\beta}\right)^{n-1}\dfrac{\left(1-\alpha^{n_1-n}\right)\left(1-\gamma^{n_2-n}\right)-\left(1-\alpha^{n_2-n}\right)\left(1-\gamma^{n_1-n}\right)}{\left(1-\beta^{n_1-n}\right)\left(1-\gamma^{n_2-n}\right)-\left(1-\beta^{n_2-n}\right)\left(1-\gamma^{n_1-n}\right)}+1 \right)\nonumber\\
	&\quad\quad+ \nu_3\left(\left(1-\beta^{n_1-n}\right)\left(1-\gamma^{n_2-n}\right)-\left(1-\beta^{n_2-n}\right)\left(1-\gamma^{n_1-n}\right)\right).
\end{align}	

As before, we first estimate
$
	\nu_p\left(\left(1-\beta^{n_1-n}\right)\left(1-\gamma^{n_2-n}\right)-\left(1-\beta^{n_2-n}\right)\left(1-\gamma^{n_1-n}\right)\right),	
$
for $p=2,3$. Specifically, 
\begin{align*}
	\nu_p\left(\left(1-\beta^{n_1-n}\right)\left(1-\gamma^{n_2-n}\right)-\left(1-\beta^{n_2-n}\right)\left(1-\gamma^{n_1-n}\right)\right)&\le \frac{(6(n-n_2)+30)\log \alpha}{\log p}\\
	&< \frac{\left(6\cdot 2\cdot 10^{15}(\log n)^2+30\right)\log \alpha}{\log p}\\
	&< 
	\begin{cases} 
		1.1\cdot 10^{16}(\log n)^2, & \text{if }~ p=2; \\
		6.7\cdot 10^{15}(\log n)^2, & \text{if } ~p=3,
	\end{cases}
\end{align*}
where we have used Lemma \ref{lem3.6t} and $n>310$. At this point, relations \eqref{3.29t} and \eqref{3.30t} become 
\begin{align}\label{3.32t}
	x_{min}&\le 
	\nu_2\left(\dfrac{C_\alpha}{C_\beta}\left(\dfrac{\alpha}{\beta}\right)^{n-1}\dfrac{\left(1-\alpha^{n_1-n}\right)\left(1-\gamma^{n_2-n}\right)-\left(1-\alpha^{n_2-n}\right)\left(1-\gamma^{n_1-n}\right)}{\left(1-\beta^{n_1-n}\right)\left(1-\gamma^{n_2-n}\right)-\left(1-\beta^{n_2-n}\right)\left(1-\gamma^{n_1-n}\right)}+1 \right)+ 1.1\cdot10^{16}\left(\log n\right)^2,
\end{align}
and 
\begin{align}\label{3.33t}	
	y_{min}&\le 
	\nu_3\left(\dfrac{C_\alpha}{C_\beta}\left(\dfrac{\alpha}{\beta}\right)^{n-1}\dfrac{\left(1-\alpha^{n_1-n}\right)\left(1-\gamma^{n_2-n}\right)-\left(1-\alpha^{n_2-n}\right)\left(1-\gamma^{n_1-n}\right)}{\left(1-\beta^{n_1-n}\right)\left(1-\gamma^{n_2-n}\right)-\left(1-\beta^{n_2-n}\right)\left(1-\gamma^{n_1-n}\right)}+1 \right)+ 6.7\cdot10^{15}\left(\log n\right)^2.
\end{align}	
Lastly, we estimate
\begin{align}\label{3.34t}
	\nu_p\left(\dfrac{C_\alpha}{C_\beta}\left(\dfrac{\alpha}{\beta}\right)^{n-1}\dfrac{\left(1-\alpha^{n_1-n}\right)\left(1-\gamma^{n_2-n}\right)-\left(1-\alpha^{n_2-n}\right)\left(1-\gamma^{n_1-n}\right)}{\left(1-\beta^{n_1-n}\right)\left(1-\gamma^{n_2-n}\right)-\left(1-\beta^{n_2-n}\right)\left(1-\gamma^{n_1-n}\right)}+1 \right),
\end{align}
with
$$
\lambda_{1}:=\dfrac{\alpha}{\beta}\qquad {\text{\rm and}}\qquad \lambda_{2}:=-\dfrac{C_\alpha}{C_\beta}\left(\dfrac{\left(1-\alpha^{n_1-n}\right)\left(1-\gamma^{n_2-n}\right)-\left(1-\alpha^{n_2-n}\right)\left(1-\gamma^{n_1-n}\right)}{\left(1-\beta^{n_1-n}\right)\left(1-\gamma^{n_2-n}\right)-\left(1-\beta^{n_2-n}\right)\left(1-\gamma^{n_1-n}\right)}\right),
$$ 
as before. We again distinguish two cases.

\medskip

{\bf Case 1.} {\it $\lambda_1$ and $\lambda_2$ are multiplicatively independent. }

\medskip

We apply Lemma \ref{lem:Bugt} to \eqref{3.34t} with the field $\mathbb{Q}(\alpha,\beta)$ of degree $D := 6$. Since $h(\lambda_{1})=3(\log \alpha)/3=\log\alpha$, we can choose $h(\lambda_{1}):=\log\alpha$ for all cases $p=2,3$. Like in Subsection \ref{subsec3.2}, Case 1 of Lemma \ref{lem3.3t},
\begin{align*}
	h'(\lambda_{2}) 
	& <  4+\frac{8(n-n_2+4)}{3}\log \alpha\\
	& <  4+\frac{8\left(2\cdot 10^{15}(\log n)^2+4\right)}{3}\log \alpha\\
	&<3.3\cdot 10^{15}(\log n)^2,
\end{align*}
where we have used Lemma \ref{lem3.6t} and $n>310$. Therefore,
\[
E'=\dfrac{b_1}{h'(\lambda_{2})}+\dfrac{b_2}{h'(\lambda_{1})}=\dfrac{n-1}{h'(\lambda_{1})}+\dfrac{1}{h'(\lambda_{2})}<n.
\]
We have
\begin{align*}
	E&=\max\left\{\log E'+\log\log p+0.4, 10, 10\log p\right\}\\
	&=\max\left\{\log n +\log\log p+0.4, 10,10\log p\right\}.
\end{align*}
If $n>36100$, then $E<1+\log n$ in both cases $p=2,3$. Therefore, Lemma \ref{lem:Bugt} gives
\begin{align*}
	&\nu_p\left(\dfrac{C_\alpha}{C_\beta}\left(\dfrac{\alpha}{\beta}\right)^{n-1}\dfrac{\left(1-\alpha^{n_1-n}\right)\left(1-\gamma^{n_2-n}\right)-\left(1-\alpha^{n_2-n}\right)\left(1-\gamma^{n_1-n}\right)}{\left(1-\beta^{n_1-n}\right)\left(1-\gamma^{n_2-n}\right)-\left(1-\beta^{n_2-n}\right)\left(1-\gamma^{n_1-n}\right)}+1 \right)\\
	&\le \dfrac{24pg}{(p-1)(\log p)^4} E^2D^4   h'(\lambda_{1})h'(\lambda_{2}  )\\
	&<\dfrac{24pg}{(p-1)(\log p)^4} (\log n)^2\left(1+\frac{1}{\log 36100}\right)^2\cdot 6^4  ( \log\alpha )\cdot 3.3\cdot10^{15}(\log n)^2\\
	&<\dfrac{7.6\cdot 10^{19}pg}{(p-1)(\log p)^4} \left(\log n\right)^4.
\end{align*}
Hence, inequalities \eqref{3.32t} and \eqref{3.33t} become 
\begin{align*}
	x_{min}&< 
	\dfrac{7.6\cdot 10^{19}\cdot 2\cdot 3}{(2-1)(\log 2)^4} \left(\log n\right)^4+ 1.1\cdot10^{16}\left(\log n\right)^2<2\cdot 10^{21} (\log n)^4,
\end{align*}
and 
\begin{align*}
	y_{min}&<\dfrac{7.6\cdot 10^{19}\cdot 3\cdot 13}{(3-1)(\log 3)^4} \left(\log n\right)^4+ 6.7\cdot10^{15}\left(\log n\right)^2<1.1\times 10^{21} (\log n)^4.
\end{align*}	
In the above, we used $g=3$ when $p=2$ and $g=13$ when $p=3$.

\medskip

{\bf Case 2.} {\it $\lambda_1$ and $\lambda_2$ are multiplicatively dependent.}

\medskip

Like the explanation done before in Subsection \ref{subsec3.2},
$$
\nu_p(\lambda_1^{n-1}\lambda_2+1)\le 1+\frac{\log n^2}{\log p}<3\log n,
$$
for $n>36100$ and $p\in \{2,3\}$. Comparing this to Case 1, we see that the bounds from Case 1 still hold. 	
This completes the proof of Lemma \ref{lem3.7t}.
\end{proof}
To conclude this subsection, we consider a fourth solution $(n_3,x_3,y_3)$, with $n>n_1>n_2>n_3$ and we find an absolute bound for $n$. We prove the following result.
\begin{lemma}\label{lem3.8t}
	If $c\le -1$ and $n>36100$,  then
	\[
	n<6.1\cdot 10^{29},\qquad x<5.5\cdot 10^{29}, \qquad y<3.5\cdot 10^{29}.
	\]
\end{lemma}
\begin{proof}
Lemma \ref{lem3.7t} tells that out of any four solutions, the minimal $x$ and $y$ are bounded by $2\cdot  10^{21}(\log n)^4$. So, out of the five solutions, at most one of them has $x$ which is not bounded by $2\cdot 10^{21}(\log n)^4$ and at most one of them has $y$ which is not bounded by $2\cdot 10^{21}(\log n)^4$. Hence, at least one of the solutions has both $x$ and $y$ bounded by $2\cdot 10^{21}(\log n)^4$, which in particular, shows that,
	$$
	X_4=x_4\log 2+y_4\log 3<(\log 2+\log 3)\cdot 2\cdot 10^{21}(\log n)^4<3.6\cdot 10^{21} (\log n)^4.
	$$
Now, by Lemmas \ref{lem2.9t} and \ref{lem3.5t}, we get
	\begin{align*}
		n\log\alpha-3&<X=X_4+(X_3-X_4)+(X_2-X_3)+(X_1-X_2)+(X-X_1)\\
		&< 3.6\cdot10^{21}(\log n)^4 +4\cdot 5.55\cdot 10^{15} \log n\\
		&<3.61\cdot 10^{21} \left(\log n\right)^4,
	\end{align*} 	
	which implies
	\begin{align}\label{3.35t}
		\dfrac{n}{(\log n)^4}<6\cdot 10^{21}.
	\end{align}	
	We apply Lemma \ref{Guzt} to inequality \eqref{3.35t} above with $ z:=n $, $ m:=4 $, $T:=6\cdot10^{21}$.
	Since $T>(4\cdot 4^2)^4$, we get
	$$n<2^m T(\log T)^m = 2^4 \cdot 6\cdot 10^{21}(\log 6\cdot 10^{21})^4 < 6.1\cdot 10^{29}.$$	
	Further, we have by Lemma \ref{lem2.9t} that
	\begin{align*}
		X&<n\log\alpha+60\left(\log (n\log\alpha)\right)^2;\\
		x\log2+y\log3&<6.1\cdot 10^{29}\log\alpha+60\left(\log (6.1\cdot 10^{29}\log\alpha)\right)^2<3.8\cdot 10^{29}.
	\end{align*}	
	This gives $x<5.5\cdot 10^{29}$ and $y<3.5\cdot 10^{29}$ and completes the proof of Lemma \ref{lem3.8t}.
\end{proof}

\subsection{Reduction of the upper bound on $n$}
Again, we use the LLL--reduction methods, the theory of continued fractions and as well as $p-$adic reduction methods due to \cite{PET} to find a rather small bound for $n$.

First, we consider \eqref{3.25t} with the assumption that $X-X_1\ge 2$ and take logarithms, that is
\[
|\Lambda_2|=\left|\log C_\alpha +(n-1)\log \alpha-x\log 2-y\log3\right|<4.5\exp(-(X-X_1)).
\]
Like before, we consider the approximation lattice
$$ \mathcal{A}=\begin{pmatrix}
	1 & 0  & 0 \\
	0 & 1 & 0 \\
	\lfloor M\log 2\rfloor & \lfloor M\log 3\rfloor& \lfloor M\log\alpha \rfloor
\end{pmatrix} ,$$
with $C:= 10^{90}$ and choose $v:=\left(0,0,-\lfloor M\log C_\alpha \rfloor\right)$. Now, by Lemma \ref{lem2.5t}, we get  $|\Lambda_2|>c_1:=9.3\cdot 10^{-32}$ and hence $c_2:=5.12\cdot 10^{30}$. Moreover, by Lemma \ref{lem3.8t}, we have
\[
x<A_1=5.5\cdot 10^{29},~~y<A_2=3.5\cdot 10^{29},~~n-1<n<A_3=6.1\cdot 10^{29}.
\]
So, Lemma \ref{lem2.6t} gives $S=4.25\cdot 10^{59}$ and $T=7.55\cdot 10^{29}$. Since $c_2^2\ge T^2+S$, then choosing $c_3:=4.5$ and $c_4:=1$, we get $X-X_1\le 138$.

Next, we continue with the assumption that $n>n_1$, and consider the inequality
\[
0<2^x3^y-2^{x_1}3^{y_1}=T_n-T_{n_1}<T_n<\alpha^{n-1}.
\]
Dividing through by $2^x3^y$ and taking logarithms, we get
\[
|\tau|:=\left|(x-x_1)\log 2+(y-y_1)\log3\right|<1.5\dfrac{\alpha^{n-1}}{2^x3^y},
\]
where we have assumed that $2^x3^y>\alpha^n$ and applied Lemma \ref{lem2.1t}. Next, we divide the above equation by $|y_1-y|\log 2$ and get
\begin{align*}
	\left|\dfrac{\log 3}{\log 2} - \dfrac{x-x_1}{y_1-y} \right|<\dfrac{1.2\alpha^n}{2^x3^y|y_1-y|},
\end{align*}
since $\alpha>1.83$ and $y_1$, $y$ are distinct. Note that $y_1$ and $y$ are indeed distinct since if they were not, then $0.5\le |e^\tau-1|=\left|1-2^{x_1-x}\right|< \alpha^{n-1}/2^x3^y<0.5$, a contradiction. By Lemma \ref{lem:Legt} with $\mu:=\dfrac{\log 3}{\log 2} $ and $M:=10^{48}$, we have 
\begin{align*}
	\dfrac{1}{(a(M)+2)(y_1-y)^2}<\left|\dfrac{\log 3}{\log 2} - \dfrac{x-x_1}{y_1-y} \right|<\dfrac{1.2\alpha^n}{2^x3^y|y_1-y|},
\end{align*}
where $a(M)=55$ (in fact, $q_{100}>10^{48}$ and $\max\{a_k: 0\le k\le 100\}=55$). Multiplying the above inequality by $|y_1-y|\log 2$ gives
\begin{align*}
	\dfrac{\log 2}{57\cdot3.5\cdot 10^{29}}<\dfrac{\log 2}{(a(M)+2)|y_1-y|}<\left|(x-x_1)\log 2+(y-y_1)\log3\right|,
\end{align*}
so that 
\[
3.4\cdot 10^{-32}<\left|(x-x_1)\log 2+(y-y_1)\log3\right|<1.5\dfrac{\alpha^{n-1}}{2^x3^y},
\]
where we have used the upper bound $y<3.5\cdot 10^{29}$.
This gives $2^x3^y<2.5\cdot 10^{31}\alpha^n$, since $\alpha>1.83$. If the assumption that $2^x3^y>\alpha^n$ is violated, then we would have $2^x3^y\le \alpha^n<2.5\cdot 10^{31}\alpha^n$, and we are in the same situation.  Now, by assuming a third solution to \eqref{1.2t}, then
\[
0<2^x3^y-C_\alpha\alpha^{n-1}<2^{x_1}3^{y_1}+2<2.5\cdot 10^{31}\alpha^{n_1}+2<2.51\cdot 10^{31}\alpha^{n_1},
\]
which gives
\begin{align*}
	\left|2^{x} 3^{y}C_\alpha^{-1}\alpha^{-(n-1)}-1\right| <7.6\cdot 10^{31} \alpha^{-(n-n_1)}.
\end{align*}
Assume that $7.6\cdot 10^{31} \alpha^{-(n-n_1)}<0.5$, which is certainly true if $n-n_1\ge123$. Taking logarithms gives
\begin{align*}
	\left|x\log 2+y\log 3-(n-1)\log\alpha-\log C_\alpha\right|\le 1.5\cdot 10^{50}\alpha^{-(n-n_1)}.
\end{align*}
So, we consider the approximation lattice
$$ \mathcal{A}=\begin{pmatrix}
	1 & 0  & 0 \\
	0 & 1 & 0 \\
	\lfloor M\log 2\rfloor & \lfloor M\log 3\rfloor& \lfloor M\log\alpha \rfloor
\end{pmatrix} ,$$
with $M:= 10^{90}$ and choose $v:=\left(0,0,-\lfloor M\log C_\alpha \rfloor \right)$. Now, by Lemma \ref{lem2.5t}, we maintain $|\Lambda_2|>c_1:=9.3\cdot 10^{-32}$ and $c_2:=5.12\cdot 10^{30}$. Still, by Lemma \ref{lem3.8t}, we have
\[
x<A_1=5.5\cdot 10^{29},~~y<A_2=3.5\cdot 10^{29},~~n-1<n<A_3=6.1\cdot 10^{29}.
\]
So, Lemma \ref{lem2.6t} gives $S=4.25\cdot 10^{59}$ and $T=7.55\cdot 10^{29}$, so that choosing $c_3:=1.5\cdot 10^{50}$ and $c_4:=\log\alpha$, we get $n-n_1\le 417$.

To continue, we proceed as in Subsection \ref{subsec3.2.2} but with different upper bounds for $n$ and $n-n_1$. Notice that in this case, $d:=n-n_1\le 417$ and $n<6.1\cdot10^{29}<2^{99}$. We repeat the algorithm described in Subsection \ref{subsec3.2.2} after relation \eqref{3.19t} and obtain rather smaller bounds as $x_{min}, y_{min}\le 152 $.

Lastly, we find a smaller upper bound for $n$. If we write $c_X$ for the upper bound of $X-X_1$, then 
\begin{align*}
	X=X_2+(X_1-X_2)+(X-X_1)&<x_{min}\log 2 +y_{min}\log 3+2c_X;\\
	x\log 2 +y\log 3&<152\log 2 +152\log 3+2\cdot 138<549.
\end{align*}
Hence, $x<793$ and $y<500$. On the other hand, Lemma \ref{lem2.9t} implies that $n\log\alpha-3<X<549,$
so that $n<914$.

\subsection{Conclusion}
Concluding the proof of Theorem \ref{1.2ct}, we note that for \( n > 310 \), the bounds are \( n \leq 913 \), \( x < 793 \), and \( y < 500 \). To efficiently handle large \( T_n \) values, our SageMath 9.5 code utilized batch processing, iterating through all \( (n, x, y) \) combinations within these ranges. This approach selected \( c\le -1 \) values with at least five representations of the form \( T_n - 2^x 3^y \), aligning with the solutions in Theorem \ref{1.2ct}. The computation, performed on an 8GB RAM laptop, was completed in about 2 hours. \qed

\section*{Acknowledgments} 
The first author thanks the Eastern Africa Universities Mathematics Programme (EAUMP) for funding his doctoral studies. His time at Wits University, Johannesburg, greatly contributed to this paper, and he is thankful for the university's hospitality and supportive environment.  Both authors appreciate the CoEMaSS Grant \#2024--029--NUM from Wits University, which was crucial for this research.

\vspace{-0.1cm}
\section*{Addresses}

$ ^{1} $ School of Mathematics, Wits University, Johannesburg, South Africa and Department of Mathematics, School of Physical Sciences, College of Natural Sciences, Makerere University, Kampala, Uganda.

Email: \url{hbatte91@gmail.com}

\noindent 
$ ^{2} $ School of Mathematics, Wits University, Johannesburg, South Africa and Centro de Ciencias Matem\'aticas UNAM, Morelia, Mexico.

Email: \url{Florian.Luca@wits.ac.za}

\end{document}